\documentclass[a4paper,11pt]{article}

\usepackage{amsmath,amsthm,amssymb}
\usepackage{amsfonts}
\usepackage{graphicx}
\usepackage{graphics}
\usepackage{url}
\usepackage{mathrsfs}

\usepackage[dvips]{geometry}
\usepackage{algorithmic}
\usepackage{algorithm}
\usepackage[colorlinks=true,linkcolor=blue,citecolor=red]{hyperref}

\def\Tunit{\mathbf{\hat{T}}}%
\def\Bunit{\mathbf{\hat{B}}}%
\def\Nunit{\mathbf{\hat{N}}}%

\newtheorem{algo}{Algorithm}
\newenvironment{algthm}{\vskip2ex
\begin{algo}\rm\hrule\vskip1ex}
{\vskip1ex\nopagebreak\hrule\nopagebreak\end{algo}}

\newenvironment{linsys}[2][r]{%
\setlength{\arraycolsep}{.1111em}
\begin{array}{@{}*{#2}{rc}#1@{}}}{\end{array}}
\renewcommand\binom[2]{\left(\begin{array}{@{}c@{}}
#1\\[0.2ex]
#2
\end{array}\right)}

\newtheorem{theorem}{Theorem}[section]

\oddsidemargin=.cm
\def\keywordname{{\bfseries Keywords}}
\def\keywords#1{\par\addvspace\medskipamount{\rightskip=0pt plus1cm
\def\and{\ifhmode\unskip\nobreak\fi\ $\cdot$
}\noindent\keywordname\enspace\ignorespaces#1\par}}%
\def\subclassname{{\bfseries Mathematics Subject Classification
(2000)}\enspace}
\def\subclass#1{\par\addvspace\medskipamount{\rightskip=0pt plus1cm
\def\and{\ifhmode\unskip\nobreak\fi\ $\cdot$
}\noindent\subclassname\ignorespaces#1\par}}
\pdfpagewidth=\paperwidth
\pdfpageheight=\paperheight

\def\vp{{\rm d}}

\title{Cubic B\'{e}zier-Spline Curves: Interpolation and Maximum Curvature}
\author{HENK PIJLS\thanks{\parbox[t]{10cm}{\hbox{Korteweg-de Vries Institute for Mathematics, University of Amsterdam}\vskip0.2ex
\hbox{Science Park 105-107, 3rd floor (entrance via Nikhef), 1098 XG Amsterdam, The Netherlands}\vskip0.1ex
\hbox{Email:\,\url{h.g.j.pijls@uva.nl}}\vskip-0.5ex
\hbox{}}} \and LE PHUONG QUAN\thanks{
\parbox[t]{10cm}{\hbox{Department of Mathematics, College of Natural Sciences, Can Tho University}\vskip0.2ex
\hbox{3/2 Street, Can Tho City (900000), Vietnam}\vskip0.1ex
\hbox{Email:\,\url{lpquan@ctu.edu.vn}}}}}

\begin{document}
\maketitle

\begin{abstract}
In this paper, we propose a closed-form solution to the inverse problem in interpolation with periodic uniform B-spline curves. This solution is obtained by
modifying the one we have established to a similar problem with relaxed uniform B-spline curves. Then we use these solutions to determine the maximum curvature
of a B\'{e}zier-spline curve.
Our computational and graphical examples are presented with the aid of Maple procedures.
\subclass{41A15 \and 53A04 \and 65D17 \and 65D18 \and 68U05}
\end{abstract}

\section{Introduction}\label{sect1}

We refer the reader to the article \cite{twelve} that contains important and necessary information about B\'{e}zier curves and B\'{e}zier-spline curves. In that article, however,
there still remains an unsolved problem on how to find a closed $C^2$ curve that interpolates a data set of distinct points in $\mathbb{R}^2$ or even in $\mathbb{R}^3$.

Before presenting the solution to this problem, we recall fundamental results on these special curves.
Here we use position vectors, as convention, to express geometric properties in $\mathbb{R}^2$ and $\mathbb{R}^3$ as well as to do calculations with
vector functions.

\begin{figure}[htbp]
\centering
\begin{tabular}{@{}cc@{}}
\includegraphics[width=6.5cm]{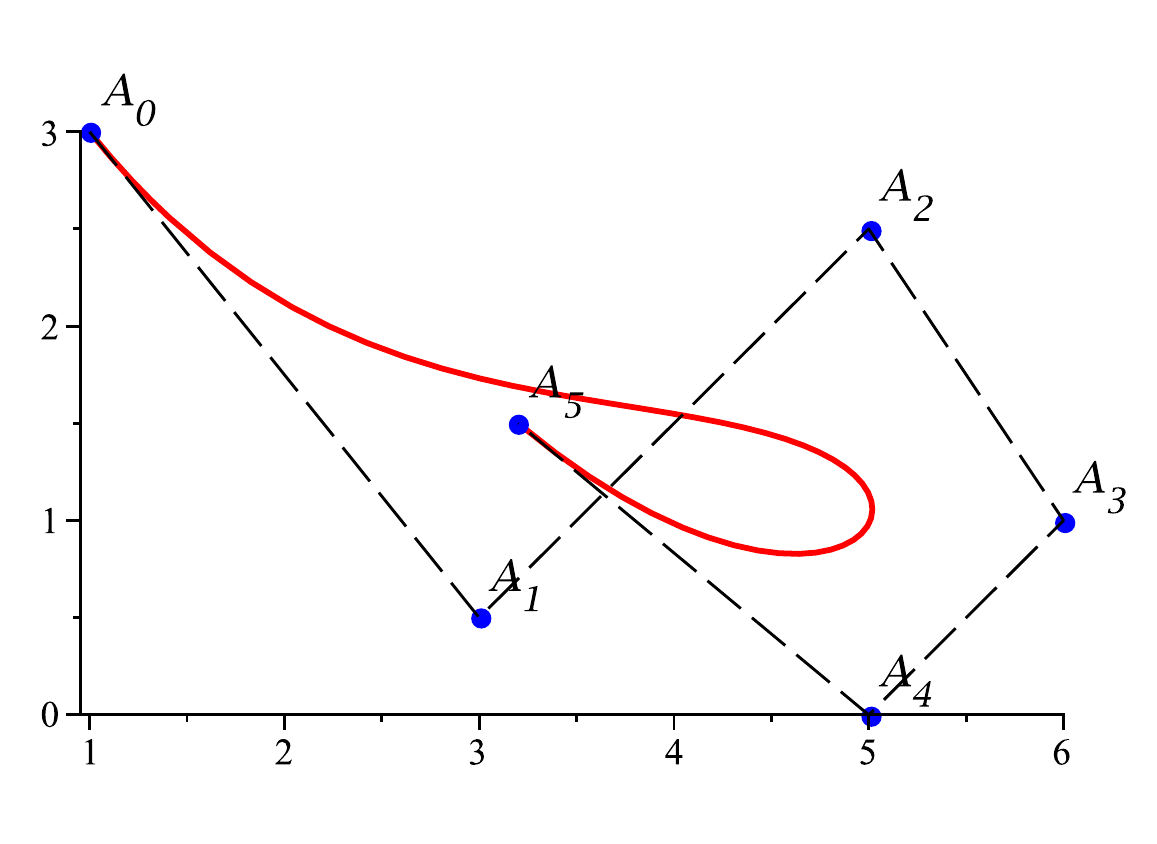}&\includegraphics[width=6.2cm]{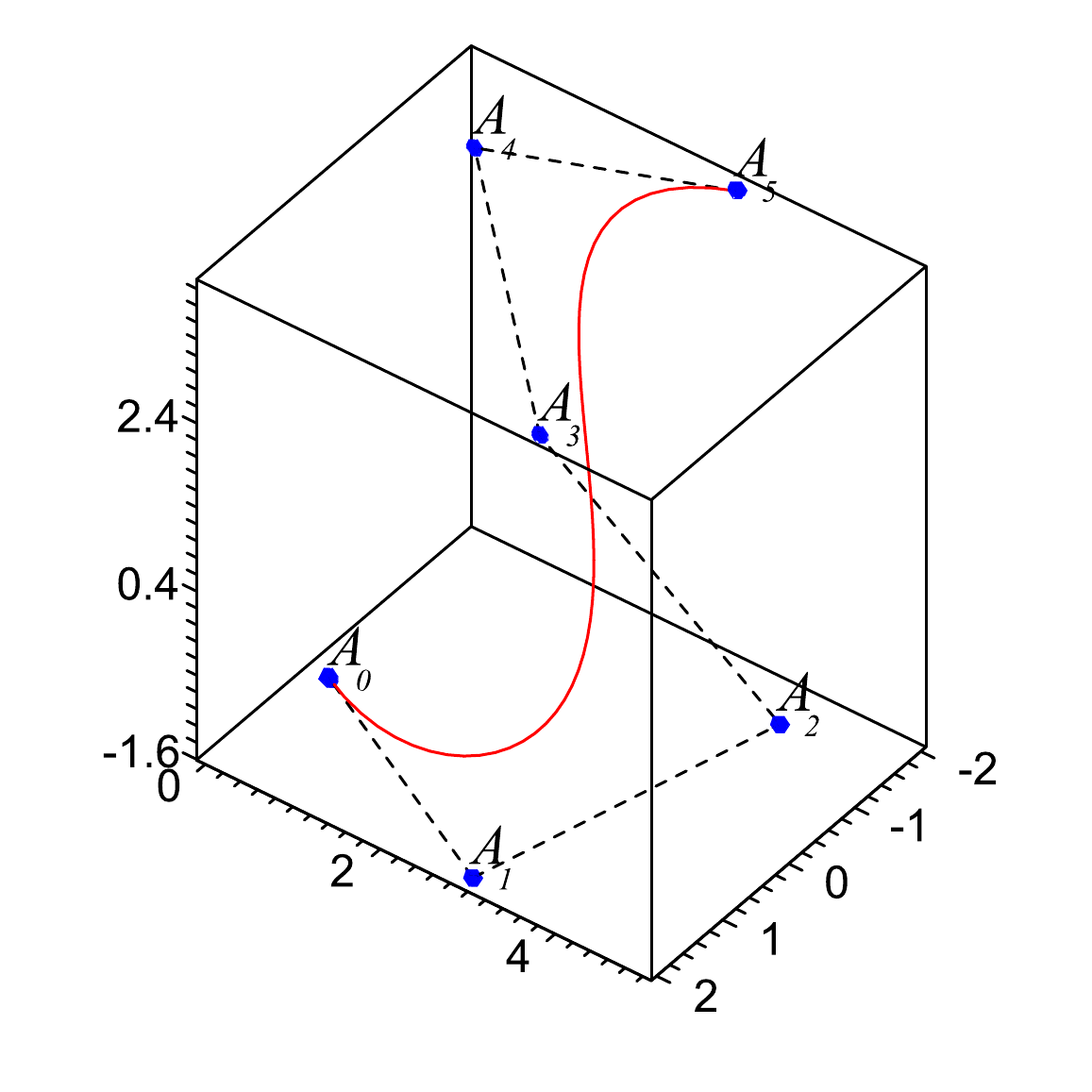}
\end{tabular}
\caption{B\'{e}zier curves of degree $5$ in $\mathbb{R}^2$ (left) and $\mathbb{R}^3$ (right).}\label{Fig1}
\end{figure}

We review a general definition of B\'{e}zier curves in $\mathbb{R}^2$ or $\mathbb{R}^3$. Let $A_0$, $A_1$, \dots, $A_n$ be $n+1$ points. A B\'{e}zier curve of degree $n\ge 2$,
taking these points as its \textsl{control points\/}, is the curve $\mathscr{C}$ that can be represented by the vector function defined on $[0,1]$:
\begin{equation}\label{eq1}
\mathbf{r}(t)=\sum_{i=0}^n\binom{n}{i}t^i(1-t)^{n-i}\mathbf{r}_i,\quad
\text{where $\binom{n}{i}=\dfrac{n!}{i!(n-i)!}$ and $0!=1$,}
\end{equation}
and $\mathbf{r}_i$ is the position vector of $A_i$, $i=0,1,\ldots,n$. We also call $\mathbf{r}(t)$ the \textsl{parametrization\/} of $\mathscr{C}$. It is easy to have that $\mathbf{r}(0)=\mathbf{r}_0$, $\mathbf{r}(1)=\mathbf{r}_n$, $\mathbf{r}'(0)=n(\mathbf{r}_1-\mathbf{r}_0)$ and
$\mathbf{r}'(1)=n(\mathbf{r}_n-\mathbf{r}_{n-1})$. All these results make the characteristics of $\mathscr{C}$ that can be described as follows: when $t$ increases from $0$ to $1$, $\mathscr{C}$
begins at point $A_0$ tangent to segment $A_0A_1$ at $A_0$, and ends at $A_n$ tangent to segment $A_{n-1}A_n$ at $A_n$. Moreover, since $0\le t\le 1$ and
$$1=[t+(1-t)]^n=\sum_{i=0}^n\binom{n}{i}t^i(1-t)^{n-i},$$
$\mathscr{C}$ is by definition contained within the convex hull of all the control points. In the case of $n = 3$, $\mathscr{C}$ is called a \textsl{cubic B\'{e}zier curve\/} and
we derive its vector function from the expanded form of (\ref{eq1}) as
\begin{equation}\label{eq2}
\mathbf{r}(t)=(1-t)^3\mathbf{r}_0+3t(1-t)^2\mathbf{r}_1+3t^2(1-t)\mathbf{r}_2+t^3\mathbf{r}_3.
\end{equation}
In Figure \ref{Fig1}, we give B\'{e}zier curves in $\mathbb{R}^2$ and $\mathbb{R}^3$.

\begin{figure}[htbp]
\centering\includegraphics[width=6.5cm]{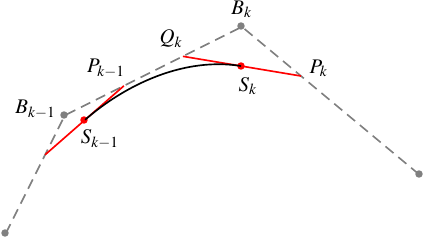}
\caption{The part $\mathscr{C}_k$ of $\mathscr{C}$ with its control points $S_{k-1}$, $P_{k-1}$, $Q_k$, $S_k$.}\label{Fig2}
\end{figure}

Another kind of B\'{e}zier curves called a \textsl{cubic and uniform B\'{e}zier-spline\/} curve with ordered control points $B_0$, $B_1$, \dots, $B_n$ ($n\ge 3$), beginning at $B_0$ and ending
at $B_n$, can be described as follows:
\begin{itemize}
\item For $0<k<n$, divide each segment $B_{k-1}B_k$ into three equal parts with subdivision points $P_{k-1}$, $Q_k$ such that, in the direction from $B_0$ to $B_n$, each $B_k$ has $Q_k$
and $P_k$ as its immediate neighbor to the left and to the right, respectively; denote by $S_k$ the midpoint of segment $Q_kP_k$ and set $S_0=B_0$, $S_n=B_n$.
\item For $1\le k\le n$, take a cubic B\'{e}zier curve $\mathscr{C}_k$ with control points $S_{k-1}$, $P_{k-1}$, $Q_k$, $S_k$ that is represented by a formula like the one in (\ref{eq2}).
All these $\mathscr{C}_k$ are then joined to form a curve $\mathscr{C}$ that we just call here \textsl{B\'{e}zier-spline curve\/}.
A sample of $\mathscr{C}$ is given in Figure \ref{Fig2}.
\end{itemize}

\begin{figure}[htbp]
\centering
\begin{tabular}{@{}c@{}c@{}}
\includegraphics[width=7cm]{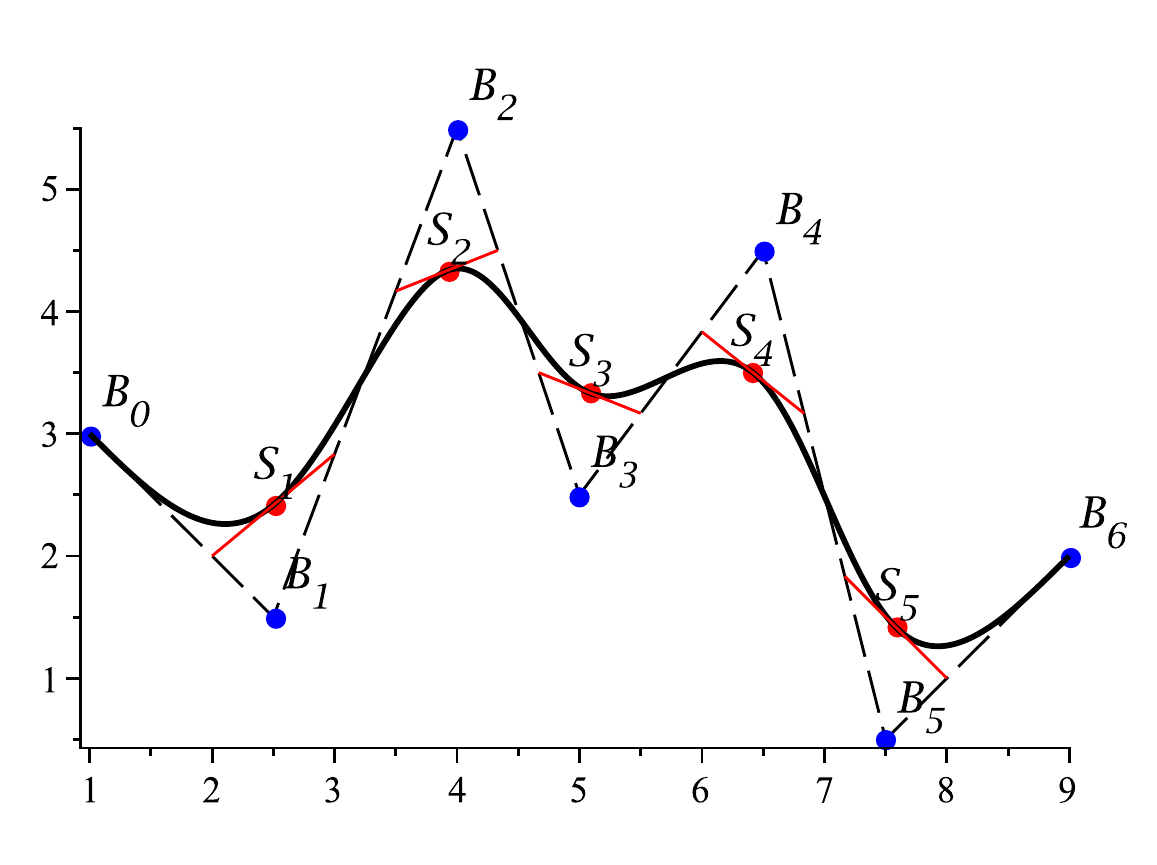}&\includegraphics[width=6cm]{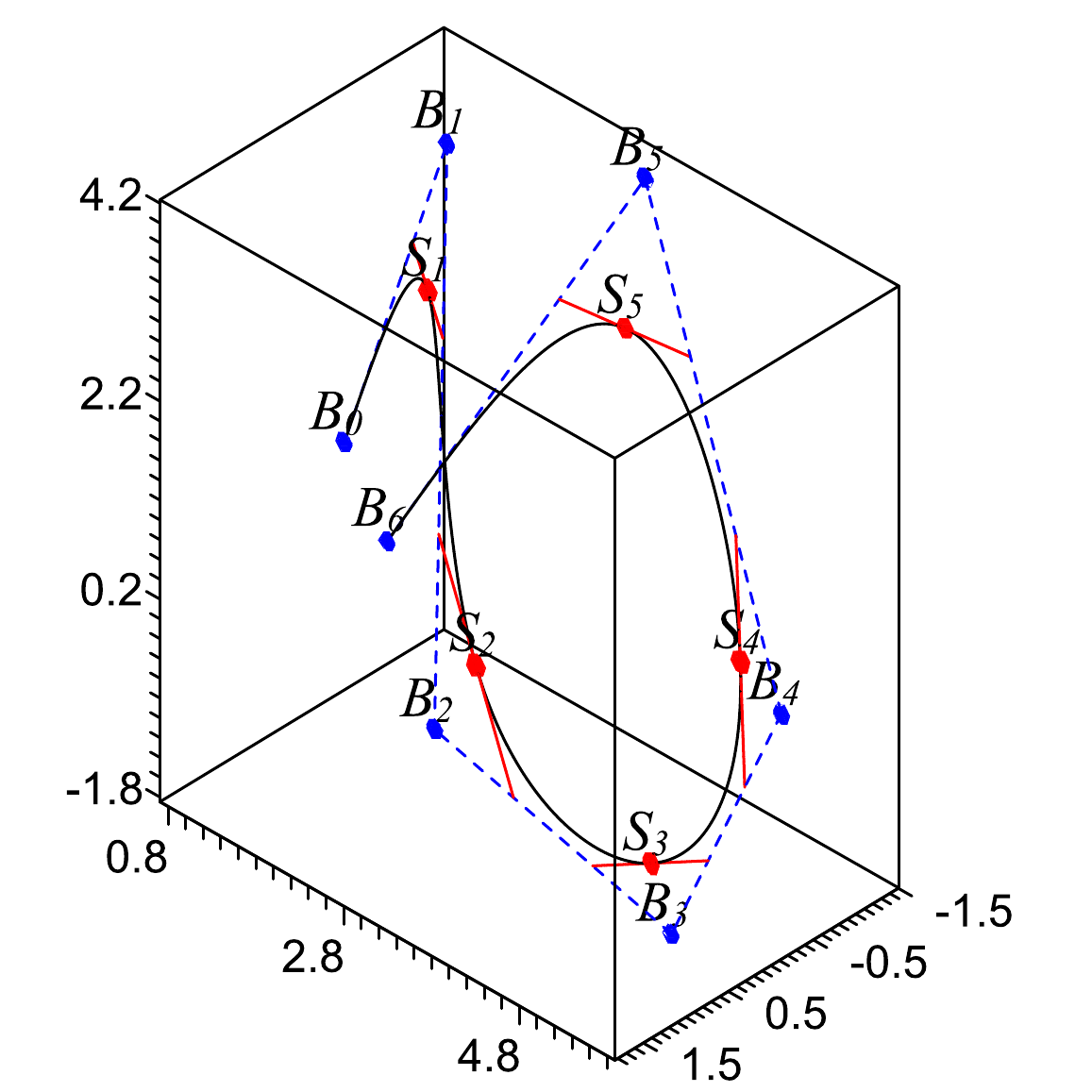}
\end{tabular}
\caption{B\'{e}zier-spline curves in $\mathbb{R}^2$ (left) and $\mathbb{R}^3$ (right).}\label{Fig3}
\end{figure}

We will give the related expressions from the above construction of a B\'{e}zier-spline curve $\mathscr{C}$. To do so, let $\mathbf{b}_0$, $\mathbf{b}_1$, \dots, $\mathbf{b}_n$
be the position vectors corresponding to $B_0$, $B_1$, \dots, $B_n$, let $\mathbf{p}_k$ be the position vector of $P_k$, $k=0,\ldots,n-1$, and let
$\mathbf{q}_k$ be the position vector of $Q_k$, $k=1,\ldots,n$, and finally, let $\mathbf{s}_k$ be the position vector of $S_k$, $k=0,1,\ldots,n$. Then, we have
\begin{equation}\label{eq3}
\mathbf{p}_{k-1}=\mathbf{b}_{k-1}+\frac{1}{3}(\mathbf{b}_k-\mathbf{b}_{k-1}),\quad
\mathbf{q}_k=\mathbf{b}_{k-1}+\frac{2}{3}(\mathbf{b}_k-\mathbf{b}_{k-1})\quad(k=1,\ldots,n),
\end{equation}
and
$$
\mathbf{s}_k=\frac{\mathbf{p}_k+\mathbf{q}_k
}{2}\quad(k=1,\ldots,n-1),\quad \mathbf{s}_0=\mathbf{b}_0,\quad \mathbf{s}_n=\mathbf{b}_n.
$$
Moreover, from (\ref{eq2}), we derive the vector function $\mathbf{f}_k(t)$ of $\mathscr{C}_k$ on $[0,1]$ as
$$
\mathbf{f}_k(t)=(1-t)^3\mathbf{s}_{k-1}+3t(1-t)^2\mathbf{p}_{k-1}+3t^2(1-t)\mathbf{q}_k+t^3\mathbf{s}_k,\quad k=1,\ldots,n.
$$
We also have that $\mathscr{C}$ by definition begins at $S_0=B_0$ tangent to segment $B_0B_1$ at $S_0$ and ends at $S_n=B_n$ tangent to $B_{n-1}B_n$ at $S_n$,
so all the points $S_k$, $k=0,1,\ldots,n$, are on $\mathscr{C}$. To see this, we give an example on making $\mathscr{C}$ in $\mathbb{R}^2$ and $\mathbb{R}^3$ in Figure \ref{Fig3}.

Besides, it is easy to check that
$$\mathbf{f}_k'(1)=\mathbf{f}_{k+1}'(0),\quad\mathbf{f}_k''(1)=\mathbf{f}_{k+1}''(0),\quad k=1,\ldots,n-1.$$
Therefore, all parts $\mathscr{C}_k$ of $\mathscr{C}$ have the same curvature at their common points. We redefine $\mathbf{f}_k$ now by
\begin{align}
\mathbf{f}_k(t)&=(k-t)^3\mathbf{s}_{k-1}+3(t+1-k)(k-t)^2\mathbf{p}_{k-1}\nonumber\\
&\quad\quad+3(t+1-k)^2(k-t)\mathbf{q}_k+(t+1-k)^3\mathbf{s}_k, \label{eq4}
\end{align}
where $t\in[k-1,k]$. The parametrization $\mathbf{f}(t)$ of $\mathscr{C}$ is defined by $\mathbf{f}=\mathbf{f}_k$ on $[k-1,k]$, $k=1,\ldots,n$. Then $\mathbf{f}$ is a $C^2$
function on $[0,n]$.

The present paper is organized as follows. In Section \ref{sect2}, we construct the solution to the inverse problems in interpolation with B\'{e}zier-spline curves. In Section \ref{sect3}, we present an algorithm
for finding the maximum curvature of $B$-spline curves from the results obtained in Section \ref{sect2}. Our computational and graphical examples are given in Section \ref{sect4}.
In Section \ref{sect5}, we give some concluding remarks.

\section{Interpolation with B\'{e}zier-spline curves}\label{sect2}

Firstly, we consider the following problem: given an ordered set of distinct points $S_0$, $S_1$, \dots, $S_n$ in $\mathbb{R}^2$ or $\mathbb{R}^3$, how to find a $C^2$ curve that interpolates this set?
According to the characteristics of a B\'{e}zier-spline curve in Section \ref{sect1}, such a curve could be a solution to the problem. Thus, how can we find the corresponding control
points $B_0$, $B_1$, \dots, $B_n$ from the data set $\{S_0,S_1,\ldots,S_n\}$ to determine this solution with the settings $B_0=S_0$ and $B_n=S_n$? In \cite{twelve}, we have solved this inverse problem
and we present here our results in the following theorem.
\begin{theorem}\label{theorem1}
Let $\mathscr{C}$ be the B\'{e}zier-spline curve with control points $B_0$, $B_1$, \dots, $B_n$ that interpolates the data set $\{S_0,S_1,\ldots,S_n\}$ in the same meaning as introduced in Section \ref{sect1}.
If we denote by $\mathbf{s}_k$ and $\mathbf{b}_k$ the position vector of $S_k$ and $B_k$, respectively, $k=0,1,\ldots,n$, then we have $\mathbf{b}_0=\mathbf{s}_0$,
$\mathbf{b}_n=\mathbf{s}_n$, and $\mathbf{b}_k$, $k=1,\ldots,n-1$, can be determined by the formula
$$
\mathbf{b}_k=\frac{\beta_{n-1-k}}{\beta_{n-1}}\Big[(-1)^k\mathbf{s}_0+6\sum_{j=1}^{k-1}(-1)^{k-j}\beta_{j-1}\mathbf{s}_j\Big]
+\frac{\beta_{k-1}}{\beta_{n-1}}\Big[(-1)^{n-k}\mathbf{s}_n+6\sum_{j=k}^{n-1}(-1)^{j-k}\beta_{n-1-j}\mathbf{s}_j\Big],
$$
where $\beta_{\ell}$, $\ell\ge-1$, is evaluated by
$$\beta_{\ell}=\frac{(2+\sqrt{3})^{\ell+1}-(2-\sqrt{3})^{\ell+1}}{2\sqrt{3}}=2^{\ell}\sum_{m=0}^{\lfloor \ell/2\rfloor}\binom{\ell+1}{\ell-2m}(3/4)^m.$$
\end{theorem}

Therefore, the solution to the current problem can be given in the form of a vector function $\mathbf{f}(t)$, defined on $[0,n]$ by
$\mathbf{f}=\mathbf{f}_k$ on $[k-1,k]$, $k=1,\ldots,n$, where
\begin{align}
\mathbf{f}_k(t)&=(k-t)^3\mathbf{s}_{k-1}+(t+1-k)(k-t)^2(2\mathbf{b}_{k-1}+\mathbf{b}_k)\nonumber\\
&\quad\quad+(t+1-k)^2(k-t)(\mathbf{b}_{k-1}+2\mathbf{b}_k)+(t+1-k)^3\mathbf{s}_k\label{eq5}
\end{align}
is obtained from (\ref{eq3}) and (\ref{eq4}). The $C^2$ curve $\mathscr{C}$ with parametrization $\mathbf{f}(t)$ is called a \textsl{relaxed uniform B-spline curve\/}
(see \cite{seven}). Its curvatures at $S_0$ and $S_n$ are both zero, since $\mathbf{f}_1''(0)=\mathbf{f}_n''(n)=\mathbf{0}$.

\begin{figure}[htbp]
\centering
\begin{tabular}{@{}cc@{}}
\includegraphics[width=7cm]{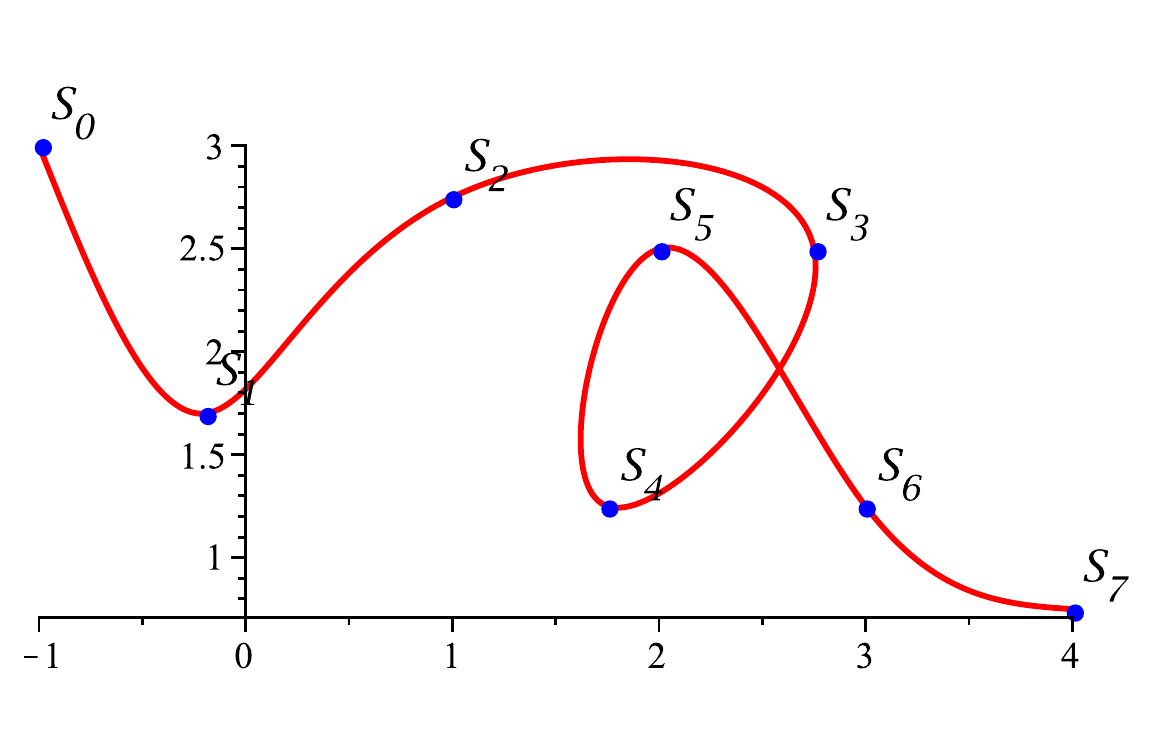}&\includegraphics[width=7.5cm]{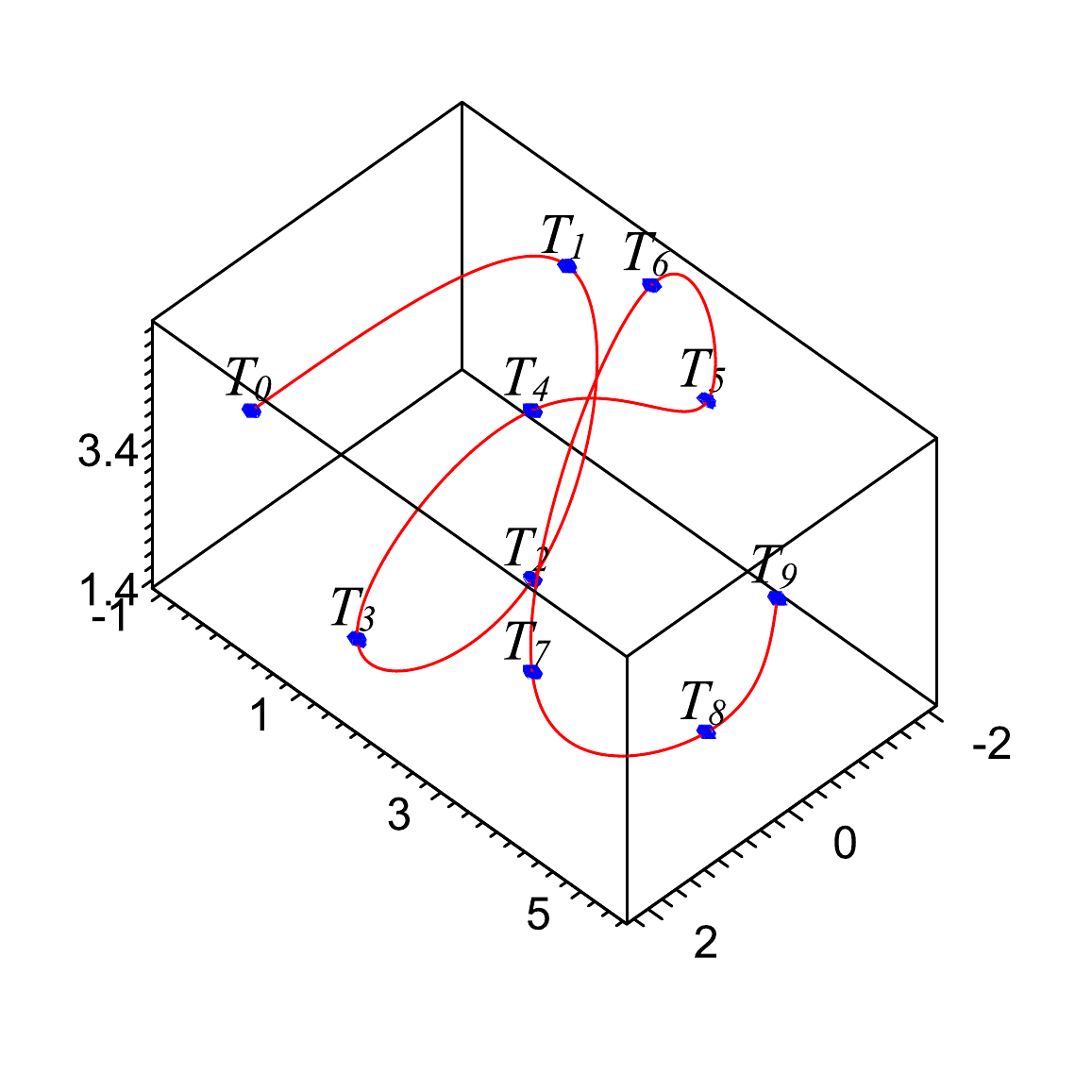}
\end{tabular}
\caption{$\mathscr{C}$ on the left, $\mathscr{L}$ on the right.}\label{Fig4}
\end{figure}

Let us give an example on interpolating data sets in $\mathbb{R}^2$ and $\mathbb{R}^3$ by relaxed uniform B-spline curves.
Choose the data sets $S$ and $T$ as follows:
\begin{align*}
S&=\{S_0=(-1,3),S_1=(-0.2,1.7),S_2=(1,2.75),S_3=(2.75,2.5),S_4=(1.75,1.25),\\
&\hskip2cm S_5=(2,2.5),S_6=(3,1.25),S_7=(4,0.75)\},\\
T&=\{T_0=(1,-1,3),T_1=(-2,0.5,4),T_2=(0,2,2),T_3=(1.5,1,1.5),T_4=(-1,1,3),\\
&\hskip1cm T_5=(-1.5,3,4.2),T_6=(-1.7,2,5),T_7=(2,4,3.5),T_8=(1,5.5,3),T_9=(-0.5,5,3.5)\}.
\end{align*}
The display of $C^2$ curves $\mathscr{C}$ and $\mathscr{L}$ that respectively interpolate $S$ and $T$ is given in Figure \ref{Fig4}.

Now, we consider a new problem: given an ordered set of distinct points $S_0$, $S_1$, \dots, $S_n$ in $\mathbb{R}^2$ or $\mathbb{R}^3$, how to
find a closed $C^2$ curve that interpolates this set?

To solve the above problem, we first follow \cite{seven} to define a closed curve. Given the control points $B_0$, $B_1$, \dots, $B_n$ of a B\'{e}zier-spline curve,
we define $B_{n+1}=B_0$ and $B_{n+2}=B_1$. Then we determine the points $P_k$, $Q_k$, $S_k$ for $k=1,\ldots,n+1$, as above, where $S_{n+1}=S_0$.
We denote by $\mathbf{b}_k$ ($k=0,\ldots,n+2$), $\mathbf{p}_k$ and $\mathbf{q}_k$ ($k=1,\ldots,n+1$), and $\mathbf{s}_k$
($k=0,\ldots,n+1$) the position vectors corresponding to the points $B_k$, $P_k$ and $Q_k$, and $S_k$; hence, we have
\begin{equation}\label{eq6}
\mathbf{p}_k=\frac{2}{3}\mathbf{b}_k+\frac{1}{3}\mathbf{b}_{k+1},\quad
\mathbf{q}_k=\frac{1}{3}\mathbf{b}_{k-1}+\frac{2}{3}\mathbf{b}_k\quad(k=1,\ldots,n+1),
\end{equation}
and
\begin{equation}\label{eq7}
\mathbf{s}_k=\frac{\mathbf{p}_k+\mathbf{q}_k
}{2}\quad(k=1,\ldots,n+1).
\end{equation}
Let $\mathscr{C}_k$ be the cubic B\'{e}zier curve with control points $S_{k-1}$, $P_{k-1}$, $Q_k$ and $S_k$, $k=1$, \dots, $n+1$, where $P_0=P_{n+1}$ so that
$$\mathbf{p}_0=\mathbf{p}_{n+1}=\frac{2}{3}\mathbf{b}_{n+1}+\frac{1}{3}\mathbf{b}_{n+2}=\frac{2}{3}\mathbf{b}_0+\frac{1}{3}\mathbf{b}_1.$$
These curves are joined together to a \textsl{periodic uniform B-spline curve\/} with control point $B_0$, \dots, $B_n$. As explained in \cite{seven}, this curve $\mathscr{C}$ is a
closed $C^2$ curve that interpolates the points $S_0$, \dots, $S_n$.

\begin{figure}[htbp]
\centering
\begin{tabular}{@{}c@{}c@{}}
\includegraphics[width=7.6cm]{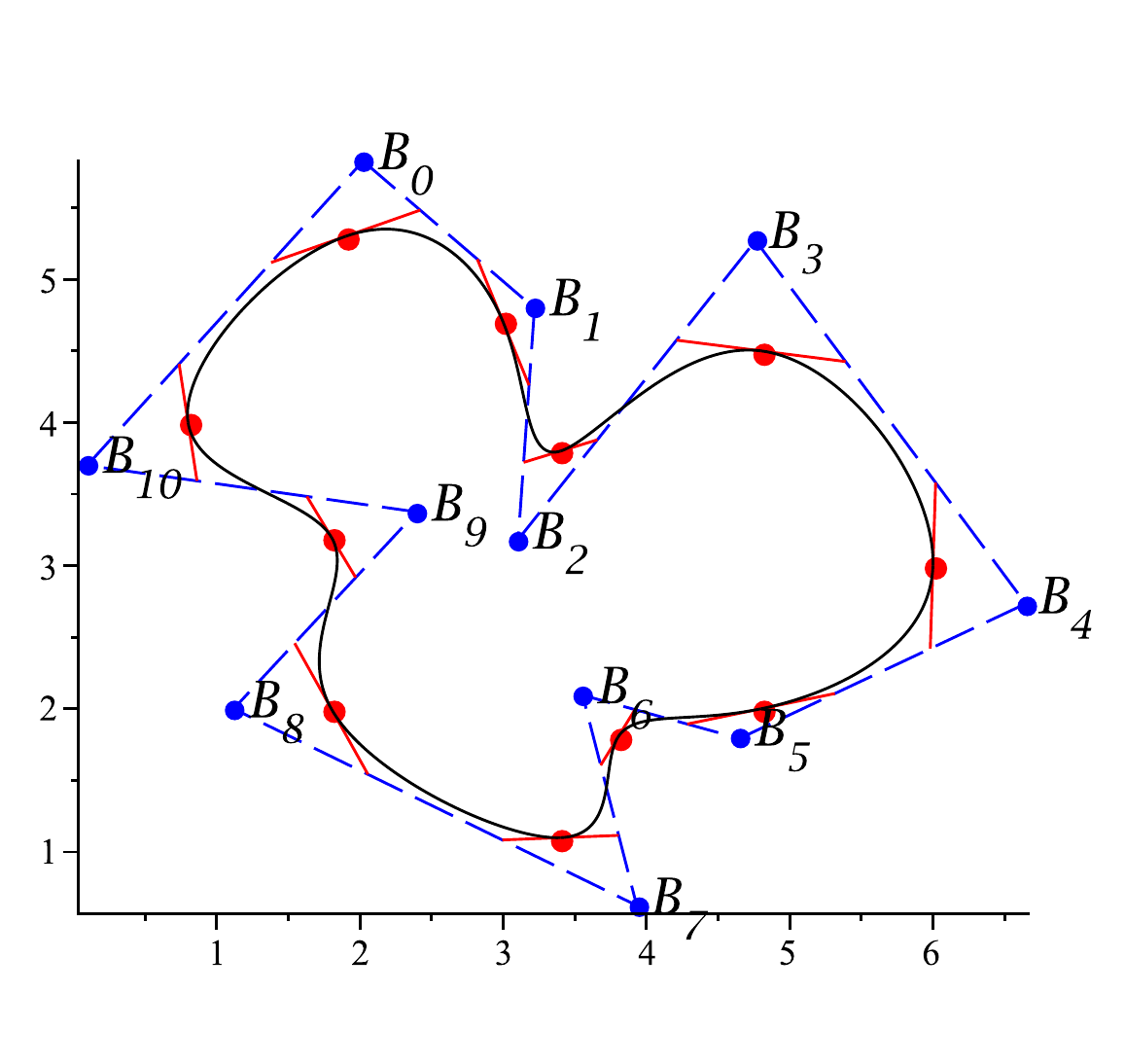}&\includegraphics[width=7.6cm]{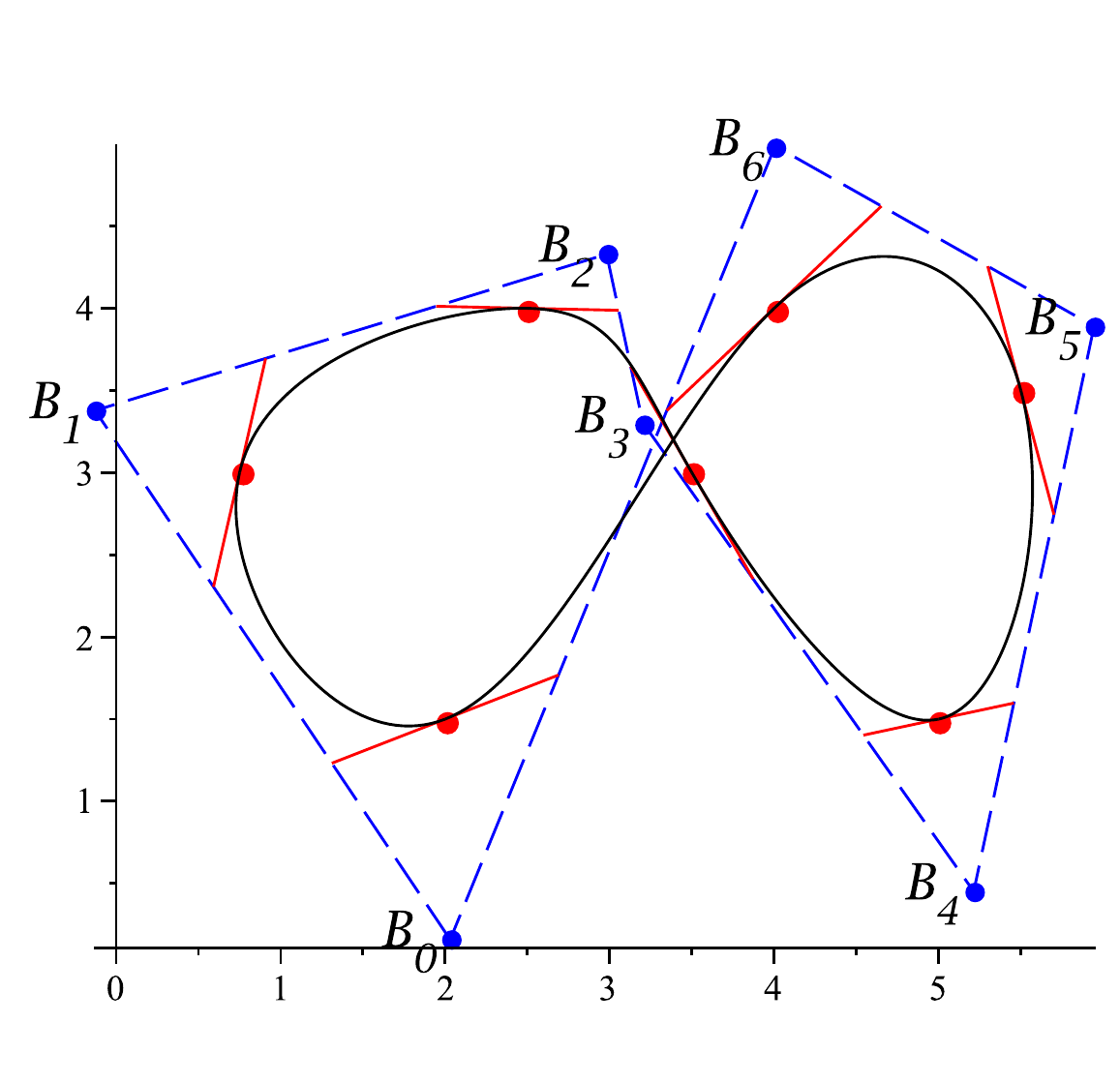}
\end{tabular}
\caption{Curves obtained from given control points $B_k$.}\label{Fig5}
\end{figure}

We give two examples of periodic uniform $B$-spline curves that interpolate sets of points $S_k$ obtained from \textsl{given\/} control points $B_k$ in $\mathbb{R}^2$.
The display of these curves is given in Figure \ref{Fig5}.

From (\ref{eq6}) and (\ref{eq7}), we derive the equations
$$\mathbf{b}_{k-1}+4\mathbf{b}_k+\mathbf{b}_{k+1}=6\mathbf{s}_k\quad(k=1,\ldots,n+1)$$
that we can write as a system of linear equations
$$\begin{linsys}[l]{6}
\mathbf{b}_0&+&4\mathbf{b}_1&+&\mathbf{b}_2&&&&&&&=&6\mathbf{s}_1\\
&&\mathbf{b}_1&+&4\mathbf{b}_2&+&\mathbf{b}_3&&&&&=&6\mathbf{s}_2\\
&&&&&&\ddots&&\ddots&&&\vdots&\\
&&&&\mathbf{b}_{n-1}&+&4\mathbf{b}_n&+&\mathbf{b}_{n+1}&&&=&6\mathbf{s}_n\\
&&&&&&\mathbf{b}_n&+&4\mathbf{b}_{n+1}&+&\mathbf{b}_{n+2}&=&6\mathbf{s}_{n+1}.
\end{linsys}$$
Since $\mathbf{b}_{n+1}=\mathbf{b}_0$, $\mathbf{b}_{n+2}=\mathbf{b}_1$ and $\mathbf{s}_{n+1}=\mathbf{s}_0$, the last equation becomes
\begin{equation}\label{eq8}
\mathbf{b}_n+4\mathbf{b}_0+\mathbf{b}_1=6\mathbf{s}_0
\end{equation}
and the remaining part of the system is now rewritten as
$$\begin{linsys}[l]{6}
4\mathbf{b}_1&+&\mathbf{b}_2&&&&&&&&&=&6\mathbf{s}_1-\mathbf{b}_0\\
\mathbf{b}_1&+&4\mathbf{b}_2&+&\mathbf{b}_3&&&&&&&=&6\mathbf{s}_2\\
&&\mathbf{b}_2&+&4\mathbf{b}_3&+&\mathbf{b}_4&&&&&=&6\mathbf{s}_3\\
&&&&&&\ddots&&\ddots&&&\vdots&\\
&&&&&&&&\mathbf{b}_{n-1}&+&4\mathbf{b}_n&=&6\mathbf{s}_n-\mathbf{b}_0
\end{linsys}$$
or
\begin{equation}\label{eq9}
\left[\begin{array}{@{}cccccc@{}}
4&1&&&&\\
1&4&1&&&\\
&1&4&1&&\\
&&\ddots&\ddots&\ddots&\\
&&&1&4&1\\
&&&&1&4
\end{array}\right]
\left[\begin{array}{@{}c@{}}
\mathbf{b}_1\\
\mathbf{b}_2\\
\mathbf{b}_3\\
\vdots\\
\mathbf{b}_{n-1}\\
\mathbf{b}_n
\end{array}\right]=
\left[\begin{array}{@{}c@{}}
6\mathbf{s}_1-\mathbf{b}_0\\
6\mathbf{s}_2\\
6\mathbf{s}_3\\
\vdots\\
6\mathbf{s}_{n-1}\\
6\mathbf{s}_n-\mathbf{b}_0
\end{array}\right].
\end{equation}
Therefore, given an ordered set of distinct points $S_0$, $S_1$, \dots, $S_n$ in $\mathbb{R}^2$ or $\mathbb{R}^3$, we need to solve (\ref{eq9}) for
$\mathbf{b}_1$, \dots, $\mathbf{b}_n$, then find $\mathbf{b}_0$ from (\ref{eq8}). In \cite{twelve}, we have solved a system like (\ref{eq9}) to obtain its unique solution
given by Theorem \ref{theorem1}. This solution can be now modified by replacing $n-1$ with $n$, $\mathbf{s}_0$ and $\mathbf{s}_n$ with $\mathbf{b}_0$, to obtain the unique solution to (\ref{eq9})
as
\begin{equation}\label{eq10}
\mathbf{b}_k=\frac{\beta_{n-k}}{\beta_n}\Big[(-1)^k\mathbf{b}_0+6\sum_{j=1}^{k-1}(-1)^{k-j}\beta_{j-1}\mathbf{s}_j\Big]
+\frac{\beta_{k-1}}{\beta_n}\Big[(-1)^{n+1-k}\mathbf{b}_0+6\sum_{j=k}^n(-1)^{j-k}\beta_{n-j}\mathbf{s}_j\Big],
\end{equation}
where $k=1$, \dots, $n$. Finally, we substitute $\mathbf{b}_1$ and  $\mathbf{b}_n$ from (\ref{eq10}) into (\ref{eq8}) and solve for $\mathbf{b}_0$.
The value of $\mathbf{b}_0$ is given in the following theorem.

{\par\noindent\textsc{Remark.\/} In our paper \cite{eleven}, we also discussed closed B\'{e}zier-spline curves.
In this article we used the equation for $\mathbf{b}_k$ given in Theorem \ref{theorem1} to obtain the formulas for
the control points of the closed curve through $n$ distinct points $S_0, \ldots, S_{n-1}$. However, we just copied this equation
without replacing $\mathbf{s}_0$ and $\mathbf{s}_n$ by $\mathbf{b}_0$. The correct formulas are now presented in the next theorem.\par}

\begin{theorem}\label{theorem2}
Let $T$ be an ordered set of distinct points $S_0$, $S_1$, \dots, $S_n$ in $\mathbb{R}^2$ or $\mathbb{R}^3$ with position vectors $\mathbf{s}_0$, \dots, $\mathbf{s}_n$. Then, we can determine
the periodic uniform B-spline curve $\mathscr{C}$ with control points $B_0$, $B_1$, \dots, $B_n$ such that their corresponding position vectors $\mathbf{b}_0$, \dots, $\mathbf{b}_n$
are defined by (\ref{eq10}), where $\mathbf{b}_0$ is given by the formula
$$
\mathbf{b}_0=\frac{3}{\beta_{n-1}-2\beta_n+(-1)^{n-1}}\Big(-\beta_n\mathbf{s}_0+\sum_{j=1}^n\big[(-1)^{j-1}\beta_{n-j}+(-1)^{n-j}\beta_{j-1}\big]\mathbf{s}_j\Big),
$$
and $\beta_{\ell}$ is given in Theorem \ref{theorem1}.
$\mathscr{C}$ interpolates the set $T$ and is composed by the cubic B\'{e}zier curves $\mathscr{C}_k$ with control points $S_{k-1}$, $P_{k-1}$, $Q_k$ and $S_k$, $k=1,\ldots,n+1$, whose parametrization $\mathbf{f}_k$ are given on $[k-1,k]$ by
\begin{equation}\label{eq11}
\mathbf{f}_k=(k-t)^3\mathbf{s}_{k-1}+3(t+1-k)(k-t)^2\mathbf{p}_{k-1}+3(t+1-k)^2(k-t)\mathbf{q}_k+(t+1-k)^3\mathbf{s}_k,
\end{equation}
where $\mathbf{p}_{k}$, $k=0,\ldots,n$, and $\mathbf{q}_k$, $k=1,\ldots,n+1$, are respectively the position vectors of $P_k$ and $Q_k$ that satisfy
$$\mathbf{p}_{k-1}=\frac{2}{3}\mathbf{b}_{k-1}+\frac{1}{3}\mathbf{b}_k,\quad \mathbf{q}_k=\frac{1}{3}\mathbf{b}_{k-1}+\frac{2}{3}\mathbf{b}_k,$$
and $\mathbf{b}_{n+1}=\mathbf{b}_0$, $\mathbf{s}_{n+1}=\mathbf{s}_0$.
\end{theorem}
If we draw the segments $P_kQ_k$, $k=1,\ldots,n+1$, we put $\mathbf{b}_{n+2}=\mathbf{b}_1$ and $\mathbf{p}_{n+1}=\mathbf{p}_0$. We also recall that our solution curve $\mathscr{C}$
has the vector function $\mathbf{f}$ that is a function of $C^2([0,n+1])$ such that $\mathbf{f}=\mathbf{f}_k$ on $[k-1,k]$, $k=1,\ldots,n+1$.

\begin{figure}[htbp]
\centering
\begin{tabular}{@{}c@{}c@{}}
\includegraphics[width=7.5cm]{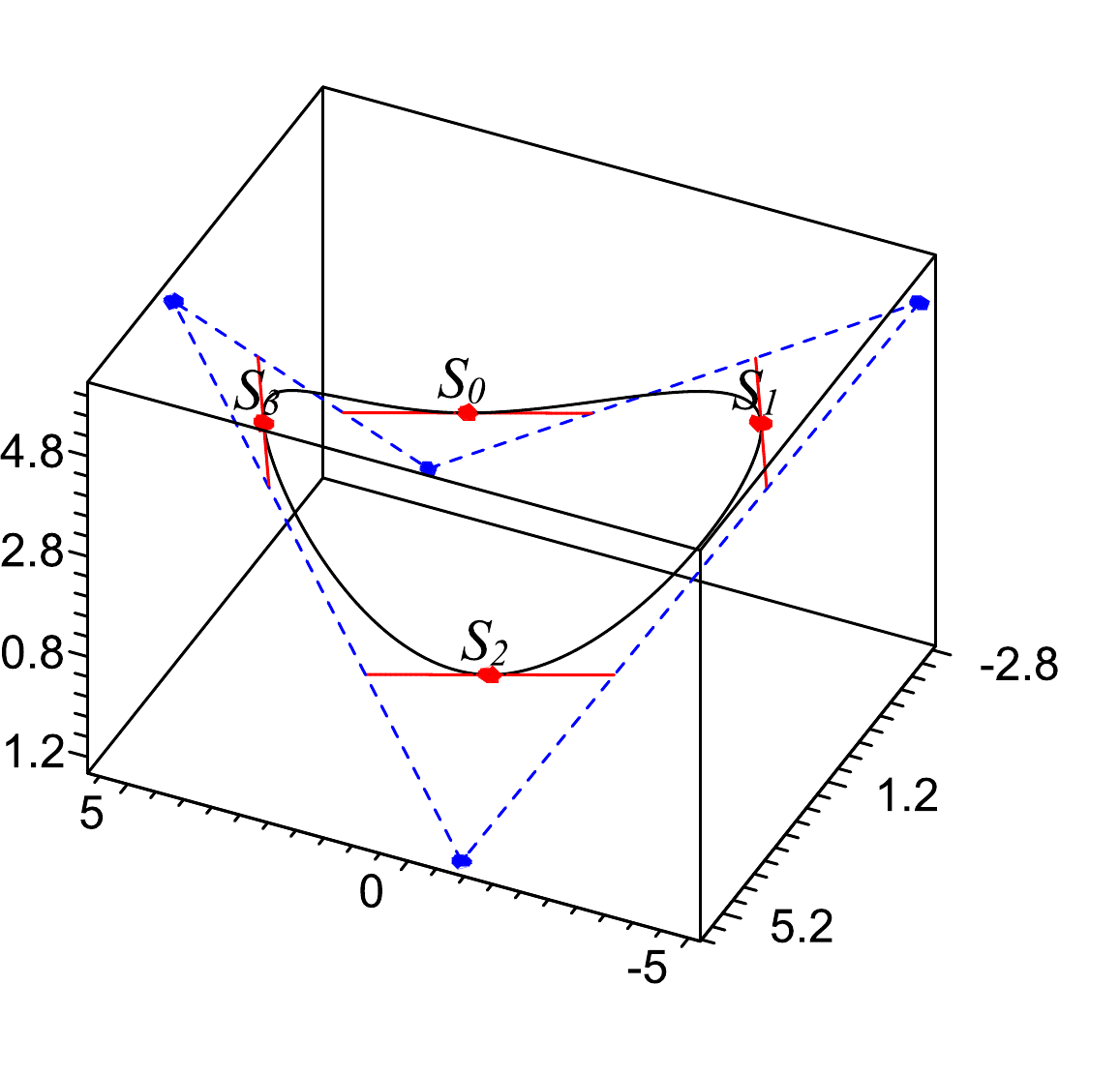}&\includegraphics[width=7cm]{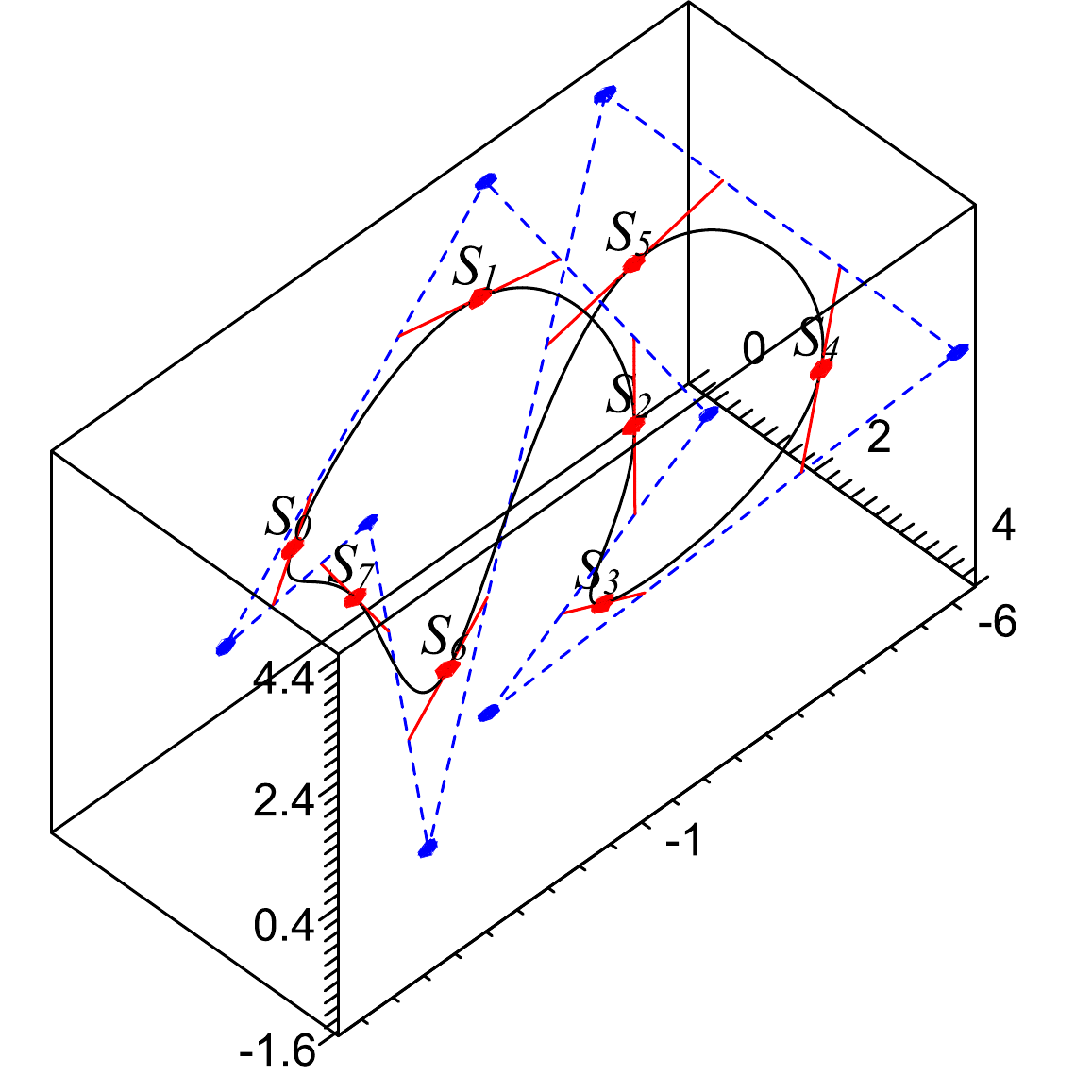}
\end{tabular}
\caption{Curves that interpolate the points $S_k$ (in red).}\label{Fig6}
\end{figure}

Finally, we give two more examples of periodic uniform $B$-spline curves that interpolate sets of \textsl{given\/} distinct points in $\mathbb{R}^3$. The display of these curves is given in Figure \ref{Fig6}.

\section{Maximum Curvature of uniform B-spline curves}\label{sect3}

In this section, we will present our algorithms for finding the maximum curvature of a B\'{e}zier-spline curve.
First we recall the definition of the \textsl{curvature\/} of a smooth curve.

At first, we consider a smooth space curve $\mathscr{L}$ parameterized by arc length $s$: $\mathbf{r}=\mathbf{r}(s)$.
The vectors
$$\mathbf{v}(s)=\mathbf{r}'(s)=\frac{\vp\mathbf{r}}{\vp s}\quad\text{and}\quad\mathbf{a}(s)=\mathbf{v}'(s)=\mathbf{r}''(s)=\frac{\vp^2\mathbf{r}}{\vp s^2}$$
are respectively called the tangent vector (or the velocity vector) and the acceleration vector of $\mathscr{L}$ at $\mathbf{r}(s)$ (or, briefly, at $s$). Because
$\mathscr{L}$ has unit speed, $\mathbf{v}(s)$ is a unit vector and is denoted by $\Tunit(s)$; hence, we have
$$\Tunit(s)=\frac{\vp\mathbf{r}}{\vp s}.$$
Then, by definition, the curvature $\kappa(s)$ of $\mathscr{L}$ at the point $\mathbf{r}(s)$ is the length of $\vp\Tunit/\vp s$. Thus,
$$\kappa(s)=\Big\|\frac{\vp\Tunit}{\vp s}\Big\|.$$

If $\kappa(s)\neq 0$, we can divide
$\vp\Tunit/\vp s$ by its length, $\kappa(s)$, to obtain a unit vector
$$\Nunit(s)=\frac{1}{\kappa(s)}\frac{\vp\Tunit}{\vp s}.$$
This vector is orthogonal to $\Tunit(s)$ and is called the unit normal to $\mathscr{L}$ at $\mathbf{r}(s)$.
At any point $\mathbf{r}(s)$ of the curve $\mathscr{L}$ where $\Tunit$ and $\Nunit$ are defined, a third unit vector, the unit binormal $\Bunit$, is defined by
$$\Bunit=\Tunit\times\Nunit.$$
The three vectors $\{\Tunit,\Nunit,\Bunit\}$ form a right-handed basis of mutually orthogonal unit vectors.
This basis is called the \textsl{Frenet frame\/} of $\mathscr{L}$ at the point $\mathbf{r}(s)$.

Now, we consider a smooth space curve $\mathscr{L}$ parameterized by a general parameter $t$: $\mathbf{r}=\mathbf{r}(t)$. Then, the unit vectors of the Frenet frame are given
by the formulas (see \cite{one})
$$\Tunit=\frac{\mathbf{v}}{\|\mathbf{v}\|},\quad\Bunit=\frac{\mathbf{v}\times \mathbf{a}}{\|\mathbf{v}\times \mathbf{a}\|},\quad \Nunit=
\Bunit\times\Tunit,$$
where $\mathbf{v}=\mathbf{v}(t)=\mathbf{r}'(t)$ and $\mathbf{a}=\mathbf{a}(t)=\mathbf{r}''(t)$. The curvature of $\mathscr{L}$ at $t$ can be now given by
\begin{equation}\label{eq12}
\kappa=\frac{\|\mathbf{v}\times\mathbf{a}\|}{v^3},
\end{equation}
where $v=\|\mathbf{v}\|$. If $\mathbf{r}(t)=(x(t),y(t),z(t))^T$, then $K=\kappa^2$ is given by
$$K=\frac{(y'z''-z'y'')^2+(z'x''-x'z'')^2+(x'y''-y'x'')^2}{[(x')^2+(y')^2+(z')^2]^3}.$$
The maximum value of $\kappa(t)$ on $[a,b]$ is the largest value at the points where $K'(t)=0$ and $t=a$, $t=b$.
However, the equation $K'(t)=0$ may be complicated.

Alternatively, if $\mathbf{r}=\mathbf{r}(s)$ is the arc length parametrization of $\mathscr{L}$, then from the relation $\kappa(s)=\Tunit'(s)\cdot\Nunit(s)$ we derive
a simple expression for $\kappa'(s)$:
$$\kappa'(s)=\Tunit''(s)\cdot\Nunit(s).$$
However, converting a general parametrization of a smooth space curve $\mathscr{L}$ into the parametrization by the arc length $s$ may be a problem.

To avoid those mentioned obstacles, many researchers restricted their attention to the class of B\'{e}zier curves and their variants.
Recently, the papers related to maximize or minimize the curvature of these curves provided a lot of
theoretical results and useful algorithms, as well as practical tools. We can list here some representative papers such as \cite{two,three,four,eight,thirteen,fourteen}.
These papers and many others paid much attention to control points and polygons. Actually, these objects have direct effects
on the shape of curves, so they have to be modified in order to obtain a curve with properties needed in design applications. However, we present an alternative. Our results can be
seen as a way of approximation by interpolation with B\'{e}zier-spline curves. Therefore, we prefer to emphasize on interpolation points. These points can easily be chosen to design
curves in $\mathbb{R}^2$ or $\mathbb{R}^3$ with
desired shapes.

Now, we go back to our main purpose: finding the maximum curvature $\kappa_{\max}$ of B\'{e}zier-spline curves. We will discuss appropriate algorithms, and perform
them with Maple. Our choice of this tool is based on the power of Maple on symbolic computation and on solving the polynomial equations of high degree.
We can refer to \cite{nine,ten} and Maple help pages for details about meaning, syntax and usage of Maple commands.

Let $\mathscr{C}$ be a relaxed uniform B-spline curve that interpolates an ordered set $T$ of distinct points $S_0$, $S_1$, \dots, $S_n$ in $\mathbb{R}^3$ with corresponding position vectors
$\mathbf{s}_0$, $\mathbf{s}_1$, \dots, $\mathbf{s}_n$. Then $\mathbf{r}(t)$, the parametrization of $\mathscr{C}$,
is a piecewise cubic function of $C^2([0,n])$, given by its components $\mathbf{f}_k(t)$ in (\ref{eq5}).

Avoiding the square root function, we derive from (\ref{eq12})
$$K:=\kappa^2=\frac{\|\mathbf{v}\times\mathbf{a}\|^2}{v^6}=\frac{(\mathbf{v}\times\mathbf{a})\cdot(\mathbf{v}\times\mathbf{a})}{(\mathbf{v}\cdot\mathbf{v})^3}.$$
We have that
$$\begin{aligned}
K'&=\frac{[(\mathbf{v}\times\mathbf{a})\cdot(\mathbf{v}\times\mathbf{a})]'(\mathbf{v}\cdot\mathbf{v})^3-[(\mathbf{v}\times\mathbf{a})\cdot(\mathbf{v}\times\mathbf{a})]3(\mathbf{v}\cdot\mathbf{v})^22(\mathbf{v}\cdot\mathbf{v}')}
{(\mathbf{v}\cdot\mathbf{v})^6}\\
&=\frac{2[(\mathbf{v}\times\mathbf{a})\cdot(\mathbf{v}\times\mathbf{a}')](\mathbf{v}\cdot\mathbf{v})-6[(\mathbf{v}\times\mathbf{a})\cdot(\mathbf{v}\times\mathbf{a})](\mathbf{v}\cdot\mathbf{a})}
{(\mathbf{v}\cdot\mathbf{v})^4}.\end{aligned}$$

Then $\kappa$ attains its maximum value $\kappa_{\rm max}$ on $[0,n]$ at solutions of the equation $K'=0$ or
\begin{equation}\label{eq13}
(\mathbf{v}\cdot\mathbf{v})[(\mathbf{v}\times\mathbf{a})\cdot(\mathbf{v}\times\mathbf{a}')]-3[(\mathbf{v}\times\mathbf{a})\cdot(\mathbf{v}\times\mathbf{a})](\mathbf{v}\cdot\mathbf{a})=0
\end{equation}
in the intervals $(i-1,i)$ or at their endpoints $i-1,i$ for $i=1,\ldots,n$. We can find $\kappa_{\rm max}$ first on each interval $[i-1,i]$ by a simple pseudo-code algorithm.
It is easy to translate the statements in such an algorithm into Maple codes or into other programming languages.

\begin{algthm}\label{HQ-algo1}
Finding the maximum curvature of relaxed uniform B-spline space curves
\end{algthm}
\begin{algorithmic}[1]
\REQUIRE An ordered set $T$ of $(n+1)$ distinct points in $\mathbb{R}^3$;
\ENSURE The maximum curvature $\kappa_{\max}$ of the relaxed uniform B-spline curve interpolating $T$;
\FOR{$i=1$ \TO $n$}
\STATE $\mathbf{v}:=\mathbf{f}_i'(t)$, $\mathbf{a}:=\mathbf{f}_i''(t)$ \quad \COMMENT{$\mathbf{f}_i$ from (\ref{eq5})};
\STATE $Q:=\text{the left-hand side of (\ref{eq13})}$;
\STATE $P:=\{t\in(i-1,i)\colon Q(t)=0\}$;
\STATE $S:=\{F_i(t)\colon t\in P\}\cup\{F_i(i-1),F_i(i)\}$ \COMMENT{$F_i:=t\mapsto
\kappa(t)=\|\mathbf{v}\times\mathbf{a}\|/v^3$, $v:=\|\mathbf{v}\|$};
\STATE $m_i:=\max S$;
\STATE $A_i:=\{t\colon F_i(t)=m_i,\,\text{where $t\in P$ or $t=i-1$ or $t=i$}\}$ ;
\ENDFOR
\STATE $\kappa_{\max}:=\max\{m_1,\ldots,m_n\}$;
\STATE $A:=\{\,\}$;
\FOR{$j=1$ \TO $n$}
\IF{$m_j=\kappa_{\max}$}
\STATE $A:=A\cup A_j$;
\ENDIF
\ENDFOR
\RETURN $\kappa_{\max}$ and $A$;
\vskip1ex
\hrule
\vskip2ex
\end{algorithmic}

The most important step of Algorithm \ref{HQ-algo1} is solving (\ref{eq13}) for $t\in[i-1,i]$, $i=1,\ldots,n$. Here we use a Maple procedure to perform this task with the
calling sequence \texttt{fsolve($Q$,$t$,$i-1$..$i$)}. Note that the degree of the polynomial $Q$ is at most $7$, and the procedure \texttt{fsolve\/} is one of the most powerful tool to numerically solve such a
polynomial or ones with even higher degrees. Beside, we do need to determine the set $A_i$ at each stage $i$. If we chose the set $M$ as a collection of all set $P$ and $\{0,1,\ldots,n\}$,
the set $\{t\in M\colon F(t)=\kappa_{\rm max}\}$ would be empty (note that $F(t)=F_i(t)$ if $t\in[i-1,i]$). After the \texttt{for\/} loop, numerical methods might give different results when we recompute some value of a function at the same argument.
We can solve this problem with the declaration of \texttt{option remember\/} for our Maple procedure.

This algorithm can be modified to use for sets of points in $\mathbb{R}^2$. Let $\mathscr{L}$ be a smooth plane curve with parametrization $\mathbf{r}=\mathbf{r}(t)$.
From the formula $\kappa(s)=\Tunit'(s)\cdot\Nunit(s)$ in terms of arc length $s$, we can write $\kappa$ in terms of $t$ as
$$\kappa=\frac{1}{v}\frac{\vp \Tunit}{\vp t}\cdot\Nunit=\frac{1}{v}\Big[\Big(\frac{1}{v}\Big)'\mathbf{v}+\frac{1}{v}\mathbf{a}\Big]\cdot\Nunit=\frac{1}{v^2}|\mathbf{a}\cdot J\Tunit|=
\frac{|\mathbf{a}\cdot\mathbf{u}|}{v^3},$$
where
$$\mathbf{v}=\mathbf{r}',\quad v=\|\mathbf{v}\|=\frac{\vp s}{\vp t},\quad\mathbf{u}=J\mathbf{v},\quad\mathbf{a}=\mathbf{v}',\quad J=\left[\begin{array}{@{}cr@{}}
0&-1\\
1&0
\end{array}\right].$$
We also consider $K:=\kappa^2=(\mathbf{a}\cdot\mathbf{u})^2/(\mathbf{v}\cdot\mathbf{v})^3$ and derive the following equation from $K'=0$:
\begin{equation}\label{eq14}
(\mathbf{v}\cdot\mathbf{v})(\mathbf{a}\cdot\mathbf{u})(\mathbf{a}'\cdot\mathbf{u})-3(\mathbf{a}\cdot\mathbf{u})^2(\mathbf{v}\cdot\mathbf{a})=0.
\end{equation}
The equation (\ref{eq14}) is equivalent to (\ref{eq15}) or (\ref{eq16}):
\begin{align}
\mathbf{a}\cdot\mathbf{u}&=0, \label{eq15}\\
(\mathbf{v}\cdot\mathbf{v})(\mathbf{a}'\cdot\mathbf{u})-3(\mathbf{a}\cdot\mathbf{u})(\mathbf{v}\cdot\mathbf{a})&=0.\label{eq16}
\end{align}
These equations play the same role as (\ref{eq13}), so we can make a $2$-dimensional version of Algorithm \ref{HQ-algo1}. The expression of $Q$ is then replaced by the
left-hand sides of (\ref{eq15}) and (\ref{eq16}), say $R$ and $W$; the set
$P$ becomes the set $\{t\in(i-1,i)\colon\text{$R(t)=0$ or $W(t)=0$}\}$ and the function $\mathbf{F}_i$ is now given by
$$\mathbf{F}_i:=t\mapsto \frac{|\mathbf{a}\cdot\mathbf{u}|}{v^3}.$$

Algorithm \ref{HQ-algo1} can be modified for finding the maximum curvature of periodic uniform B-spline curves. We only replace $n$ with $n+1$ and the components
$\mathbf{f}_i$ are now taken from (\ref{eq11}). Of course, the modification depends on whether the data sets are in $\mathbb{R}^2$ or $\mathbb{R}^3$.

On the other hand, it should be necessary to have a tool for evaluating the curvature of a uniform B-spline curve at a point in the parameter interval.
The following algorithm may be used to make such a tool.

\begin{algthm}
Evaluating curvature of relaxed uniform B-spline space curves at a point
\end{algthm}
\begin{algorithmic}[1]
\REQUIRE A set $T$ of $(n+1)$ points in $\mathbb{R}^3$, a point $p\in[0,n]$;
\ENSURE The curvature $\kappa$ at $p$ of the relaxed uniform B-spline space curve interpolating $T$;
\IF{$p=n$}
\STATE $i:=n$;
\ELSE
\STATE $i:=\lfloor p\rfloor+1$;
\ENDIF
\STATE $\mathbf{v}:=\mathbf{f}_i'(t)$, $\mathbf{a}:=\mathbf{f}_i''(t)$ \quad \COMMENT{$\mathbf{f}_i$ from (\ref{eq5})};
\STATE $F:=\|\mathbf{v}\times\mathbf{a}\|/v^3$\quad \COMMENT{$v:=\|\mathbf{v}\|$};
\RETURN $F(p)$;
\vskip1ex
\hrule
\vskip2ex
\end{algorithmic}

This simple algorithm can be used for data sets of points in $\mathbb{R}^2$ and for periodic uniform B-spline curves in $\mathbb{R}^2$ or $\mathbb{R}^3$,
with appropriate modification as mentioned above.

We will use Maple procedures to perform Algorithm \ref{HQ-algo1} and its modified versions for the given data sets of points.
These sets and the related results are supplied in the next section.

\section{Examples}\label{sect4}

If the maximum curvature $\kappa_{\rm max}$ of a B\'{e}zier-spline curve is attained at $t=\alpha\in[k-1,k]$, the position $\mathbf{f}_k(\alpha)$ will be depicted in red
in the pictures below.

Our first examples are destined for data sets of points in $\mathbb{R}^2$.
\begin{enumerate}
\item Let us choose
$$T_1:=\{(-1,3),(-0.2,1.7),(1,2.75),(2.75,2.5),(1.75,1.25),(2,2.5),(3,1.25),(4,0.75)\}.$$
The relaxed uniform $B$-spline curve $\mathscr{L}_1$ interpolating $T_1$ has its maximum curvature $\kappa_{\rm max}$ at the point
$A=\mathbf{f}_5(4.012814878)$. The position of $A$ on $\mathscr{L}_1$ is given in Figure \ref{Fig7}.
\begin{figure}[htbp]
\centering
\begin{tabular}{@{}cc@{}}
\includegraphics[width=7cm]{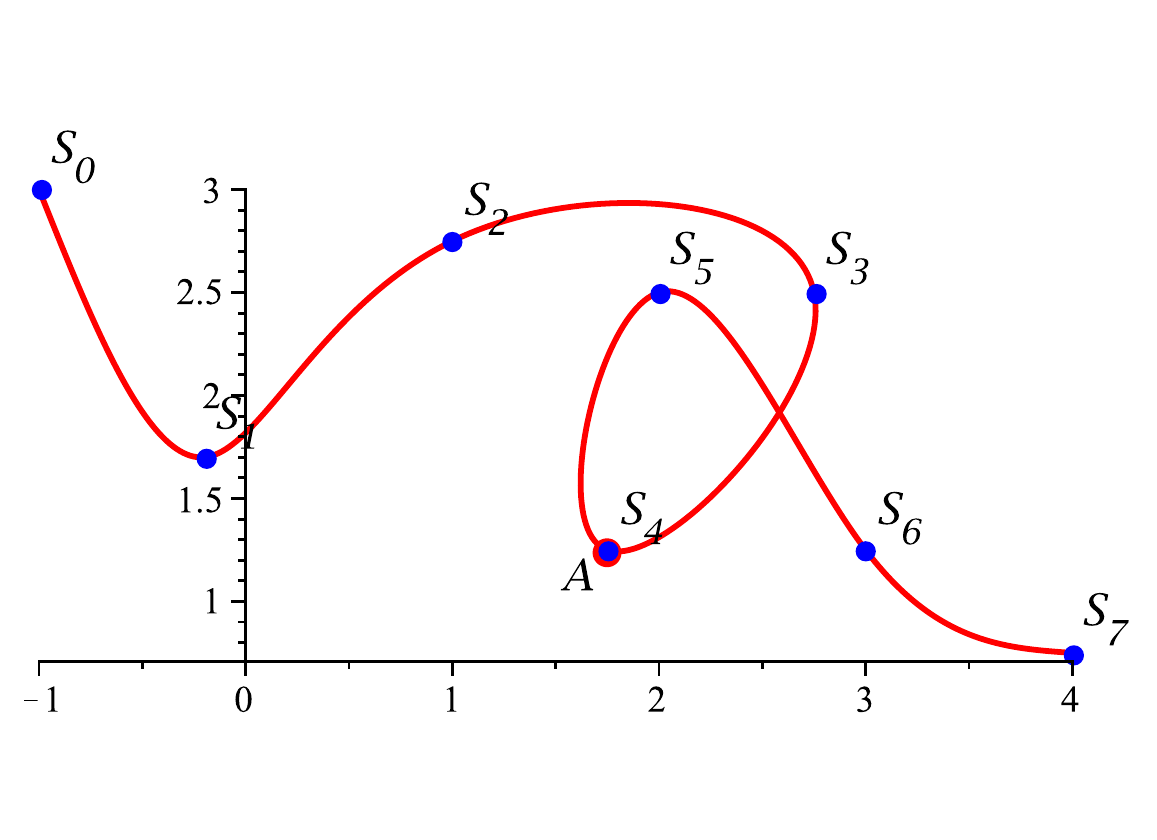}&\includegraphics[width=7cm]{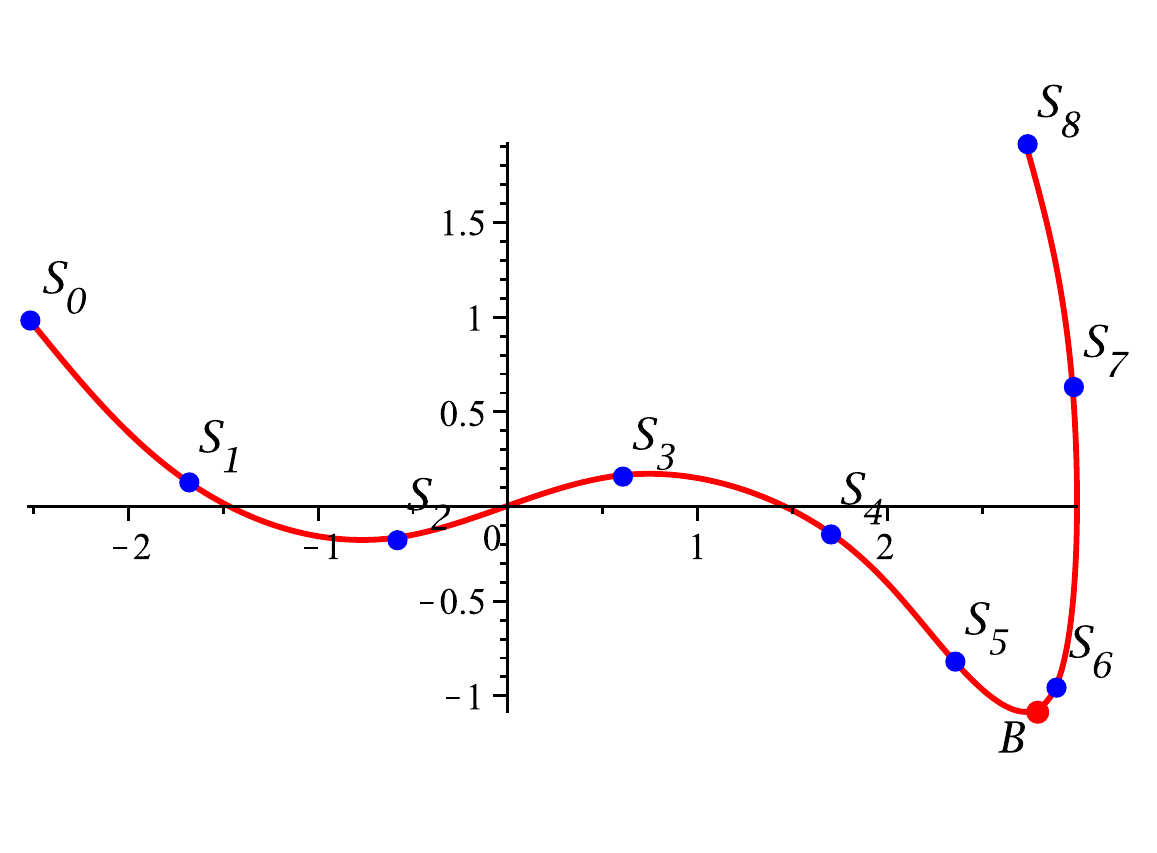}
\end{tabular}
\caption{Left: $\kappa_{\rm max}=3.263038515$ of $\mathscr{L}_1$ is attained at $A$. Right: $\kappa_{\rm max}=7.637353635$ of $\mathscr{L}_2$ is attained at $B$.}\label{Fig7}
\end{figure}
\item From the curve $\mathscr{C}$ parameterized by $\mathbf{r}(t)=(3\sin t,t\cos(3t))$ ($t\in[-1,2]$), we take
$$T_2:=\{\mathbf{r}(-1),\mathbf{r}(-0.6),\mathbf{r}(-0.2),\mathbf{r}(0.2),\mathbf{r}(0.6),\mathbf{r}(0.9),\mathbf{r}(1.3),\mathbf{r}(1.7),\mathbf{r}(2)\}.$$
The relaxed uniform $B$-spline curve $\mathscr{L}_2$ interpolating $T_2$ has its $\kappa_{\rm max}$ at the point $B=\mathbf{f}_6(5.720768304)$. The position of $B$
on $\mathscr{L}_2$ is also given in Figure \ref{Fig7}.
\item Let us take $T_3:=\{(2,1.5),(0.75,3),(2.5,4),(3.5,3),(5,1.5),(5.5,3.5),(4,4)\}$. The periodic uniform $B$-spline curve $\mathscr{L}_3$
interpolating $T_3$ has its $\kappa_{\rm max}$ at the point $S_4=\mathbf{f}_4(4)$. The position of this interpolated point
on $\mathscr{L}_3$ is shown in Figure \ref{Fig8}.
\item We choose a special set with a symmetry axis. That is the set
$$T_4:=\{(1,4),(0.6,2),(2,0.4),(3.4,1),(2.6,2.8),(2.2,2.4),(4,1.6),(4.6,3),(3,4.4)\}.$$
If $\mathscr{L}_4$ is the periodic uniform $B$-spline curve interpolating $T_4$, then it has $\kappa_{\rm max}$ at the point $C=\mathbf{f}_5(4.5000492)$. We will see this position on
$\mathscr{L}_4$ also in Figure \ref{Fig8}.
\begin{figure}[htbp]
\centering
\begin{tabular}{@{}cc@{}}
\includegraphics[width=7.5cm]{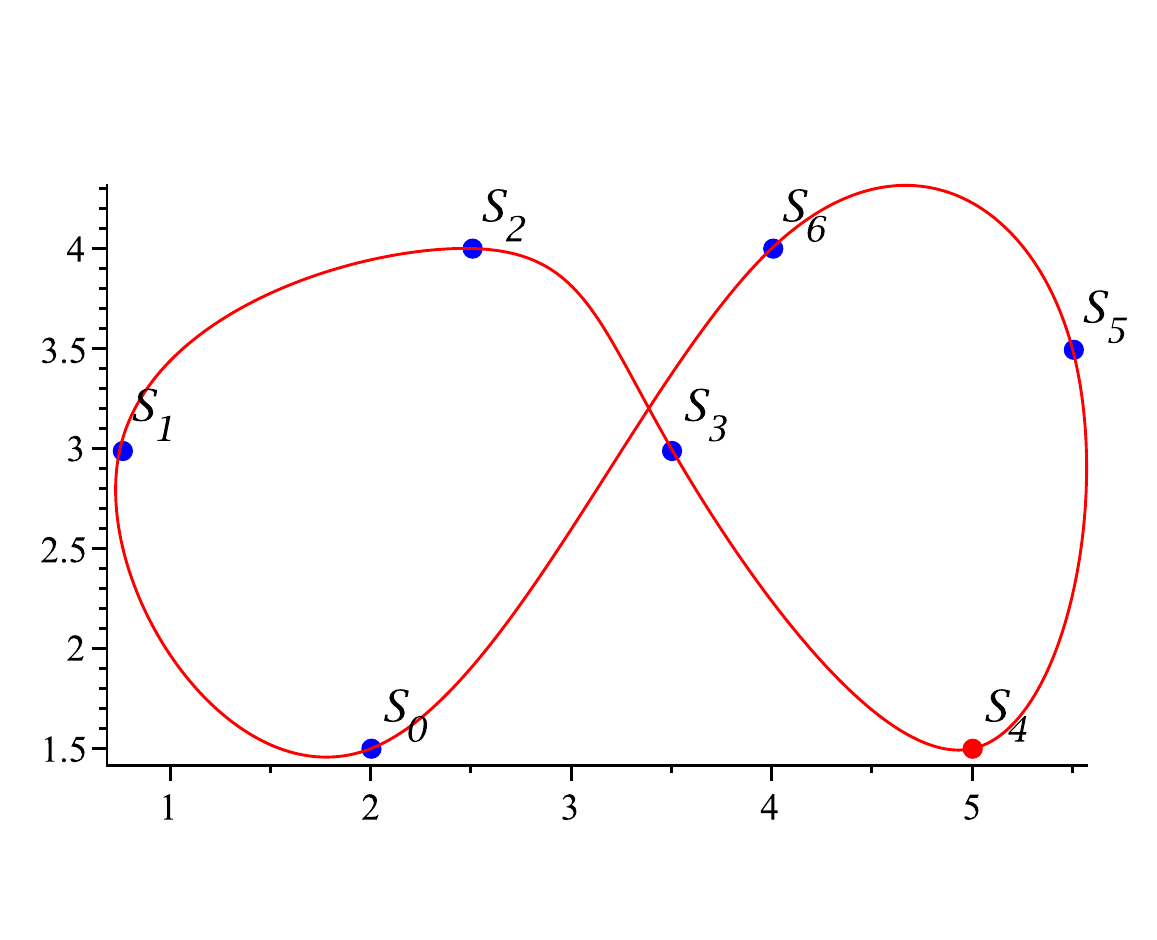}&\includegraphics[width=6.8cm]{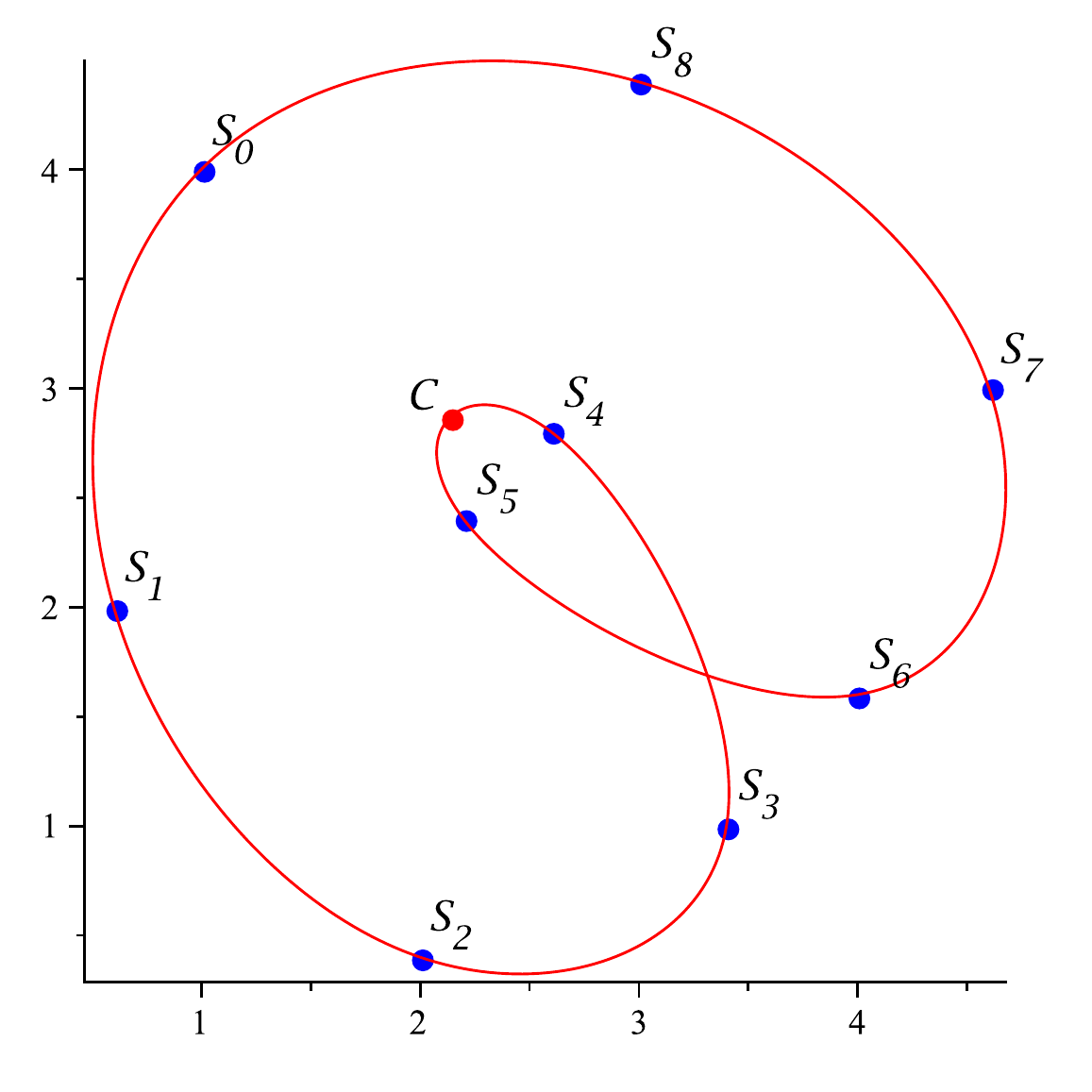}
\end{tabular}
\caption{Left: $\kappa_{\rm max}=7.98849404$ of $\mathscr{L}_3$ is attained at $S_4$. Right: $\kappa_{\rm max}=5.342692666$ of $\mathscr{L}_4$ is attained at $C$.}\label{Fig8}
\end{figure}
\item Let $\mathscr{C}$ be an ellipse parameterized by $\mathbf{r}(t)=(3\cos t,2\sin t)$, $t\in[0,2\pi]$. We take special values $t_i$ of $t$ such that the points
$\mathbf{r}(t_i)$ are symmetric about the $x$-axis. We take
\begin{align*}
T_5:&=\Big\{\mathbf{r}(0),\mathbf{r}\Big(\frac{\pi}{8}\Big),\mathbf{r}\Big(\frac{\pi}{4}\Big),\mathbf{r}\Big(\frac{3\pi}{8}\Big),\mathbf{r}\Big(\frac{5\pi}{8}\Big),
\mathbf{r}\Big(\frac{3\pi}{4}\Big),\mathbf{r}\Big(\frac{7\pi}{8}\Big),\mathbf{r}(\pi),\mathbf{r}\Big(\frac{9\pi}{8}\Big),\mathbf{r}\Big(\frac{5\pi}{4}\Big),\\
&\quad\mathbf{r}\Big(\frac{11\pi}{8}\Big),\mathbf{r}\Big(\frac{13\pi}{8}\Big),\mathbf{r}\Big(\frac{7\pi}{4}\Big),\mathbf{r}\Big(\frac{15\pi}{8}\Big)\Big\}.
\end{align*}
However, these points are not exactly symmetric about the $x$-axis, because we can only evaluate their coordinates approximately.
Then we find the maximum curvature $\kappa_{\rm max}$ of the periodic uniform $B$-spline curve $\mathscr{L}_5$ interpolating $T_5$. This value is attained at
the point $D=\mathbf{f}_{12}(11.55685766)$. The curves $\mathscr{C}$, $\mathscr{L}_5$ and the point $D$ are given in Figure \ref{Fig9}.
\begin{figure}[htbp]
\centering
\begin{tabular}{@{}c@{}}
\includegraphics[width=8cm]{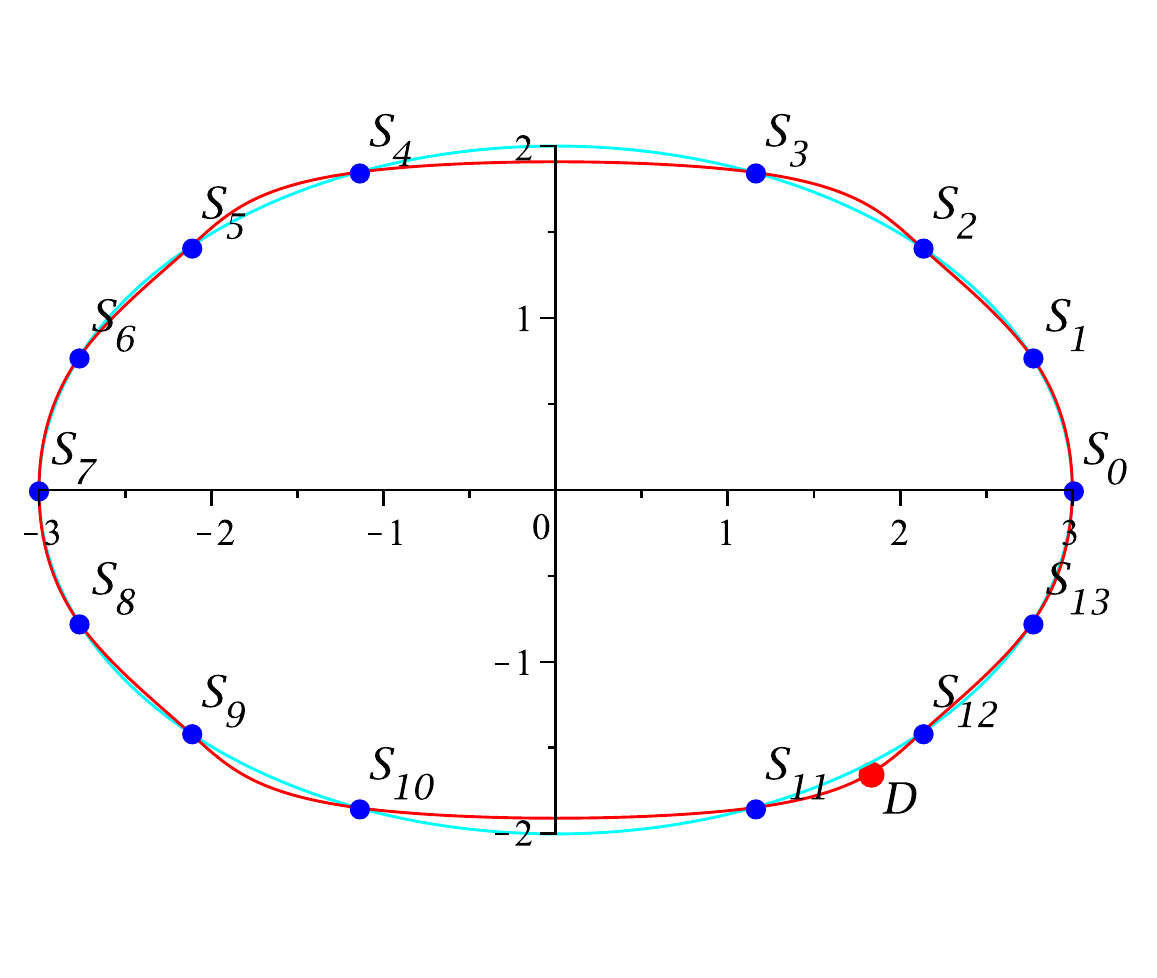}
\end{tabular}
\caption{$\mathscr{C}$ is the cyan curve. $\kappa_{\rm max}=0.95866364048$ of $\mathscr{L}_5$ is attained at $D$.}\label{Fig9}
\end{figure}
\end{enumerate}

Finally, we give examples of data sets of points in $\mathbb{R}^3$.
\begin{enumerate}
\item We first take the set
\begin{align*}
E_1:&=\{(1,2,2),(-0.5,1.5,2.5),(1,3.5,0.5),(0.5,5,-1),(-0.3,5.25,0.75),\\
&\quad(-0.75,3.5,3),(0.75,2.25,1)\}.
\end{align*}
The relaxed union $B$-spline curve $\mathscr{C}_1$ interpolating $E_1$ has its maximum curvature $\kappa_{\rm max}$ at the point $A=\mathbf{f}_1(0.8519812587)$.
The curve $\mathscr{C}_1$ and the point $A$ are given in Figure \ref{Fig10}.
\begin{figure}[htbp]
\centering
\begin{tabular}{@{}cc@{}}
\includegraphics[width=7cm]{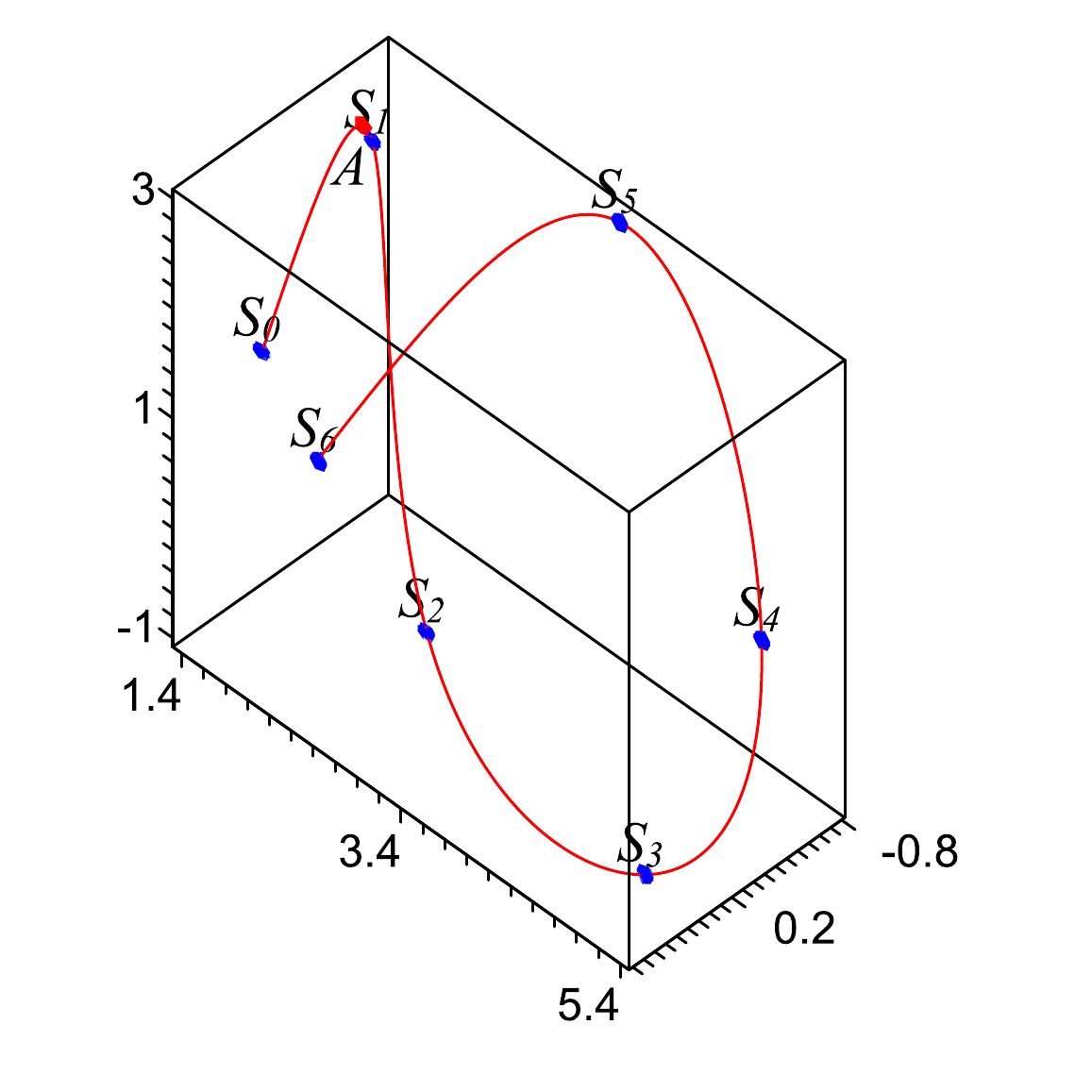}&\includegraphics[width=7cm]{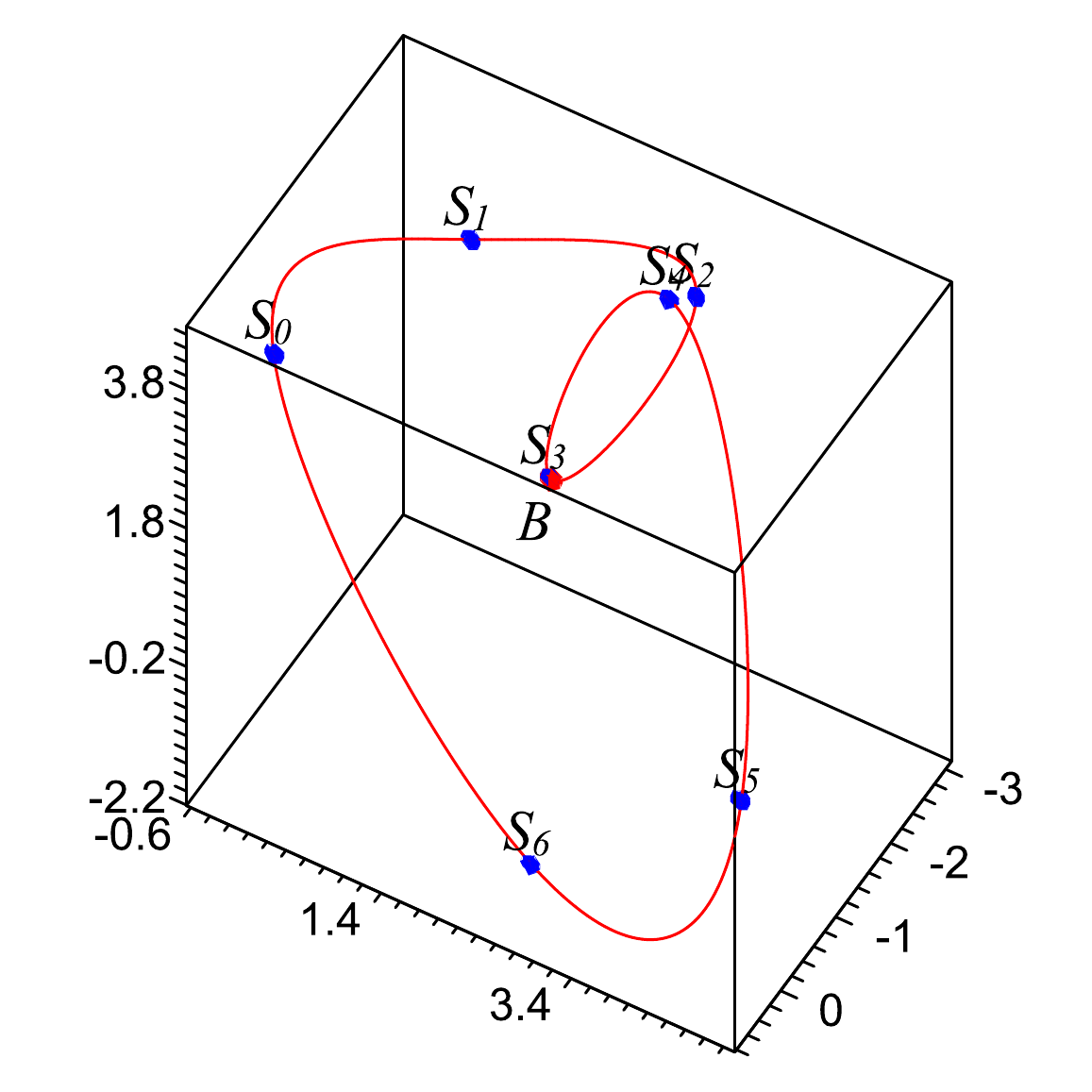}
\end{tabular}
\caption{Left: $\kappa_{\rm max}=21.16170229$ of $\mathscr{C}_1$ is attained at $A$. Right:
$\kappa_{\rm max}=10.72508237$ of $\mathscr{C}_2$ is attained at $B$.}\label{Fig10}
\end{figure}
\item Let us consider $\mathscr{C}_2$, the periodic uniform $B$-spline curve, that interpolates the data set
\begin{align*}
E_2:&=\{(-0.5,-0.5,3),(-1.5,1,4.5),(-3,2.5,3),(-1.2,2,2),(-2.5,2.5,3.5),\\
&\quad (0.5,5,1),(0,2.5,-2)\}.
\end{align*}
Then $\mathscr{C}_2$ has its $\kappa_{\rm max}$ at the point $B=\mathbf{f}_3(2.972072509)$. The curve $\mathscr{C}_2$ and the point $B$ are also shown in Figure \ref{Fig10}.
\item For the last example, we consider the closed curve $\mathscr{L}$ parameterized by
$$\mathbf{r}(t)=\big([1+0.3\cos(3t)]\cos(2t),[1+0.3\cos(3t)]\sin(2t),0.35\sin(3t))\quad(t\in[0,2\pi]\big).$$
This curve is called a \textsl{trefoil knot\/} (see \cite[p.897]{one}). We choose a set $E_3$ of points on $\mathscr{L}$ as follows
$$E_3:=\{\mathbf{r}(0),\mathbf{r}(0.4),\mathbf{r}(1.0),\mathbf{r}(1.5),\mathbf{r}(2.0),\mathbf{r}(2.5),\mathbf{r}(3.2),\mathbf{r}(3.9),\mathbf{r}(4.5),\mathbf{r}(5.1),\mathbf{r}(5.8)\}.$$
The periodic uniform $B$-spline curve $\mathscr{C}_3$ interpolating $E_3$ has its $\kappa_{\rm max}$ at the point $C=\mathbf{f}_7(6.926629821)$. The curves
$\mathscr{L}$, $\mathscr{C}_3$ and the point $C$ are given in Figure \ref{Fig11}.
\begin{figure}[htbp]
\centering
\begin{tabular}{@{}c@{}}
\includegraphics[width=8.5cm]{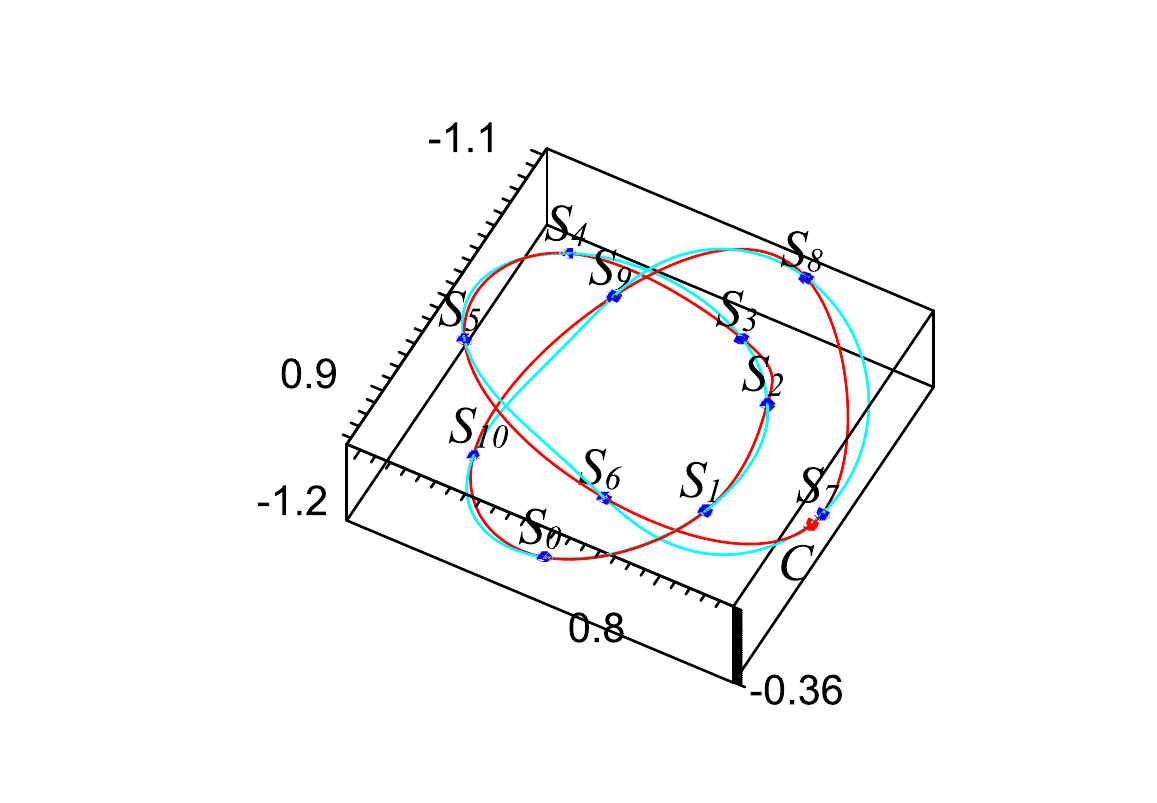}
\end{tabular}
\caption{$\mathscr{L}$ (the trefoil knot) is the cyan curve; $\kappa_{\rm max}=2.184629218$ of $\mathscr{C}_3$ is attained at $C$.}\label{Fig11}
\end{figure}
\end{enumerate}

\section{Conclusion}\label{sect5}

Our obtained results might get attention for finding curve fitting in statistics or simulating motion orbits in mechanics and computer graphics.
They could be used effectively for doing calculations to establish conditions for non-singularity of tubular surfaces and offset curves associated to a B\'{e}zier-spline curve, according to
the issues profoundly presented in \cite{five,six}. These results also give a tool to approximate the curvature extremum of a smooth curve with given parametrization.

Our interpolation vector functions can be extended directly to higher-dimensional Euclidean spaces.
Moreover, our closed-form solution of the linear system could also be applied to solve numerically differential equations using finite difference schemes.

\section*{Author Contributions}
H.~P. constructed the parametrization of B\'{e}zier-spline curves. L.~P.~Q. proposed the algorithms and wrote the Maple codes.

\section*{Acknowledgements}
The authors would like to thank the Maplesoft experts for their great work in developing Maple, a powerful and user-friendly product.

\end{document}